\documentclass[twoside,12pt]{article}
\usepackage{amsmath, epsfig}
\newcommand{\R}{\mathbf  R}

\numberwithin{equation}{section}

\newcommand{\Alphatotal}[0]{\overline{A}}
\newcommand{\Omegatotal}[0]{\overline{\Omega}}
\newcommand{\be}{\begin{equation}}
\newcommand{\ee}{\end{equation}}
\newcommand{\f}{\frac}
\newtheorem{theorem}{Theorem}[section]
\newcommand{\lb}{\lambda}
\newcommand{\p}{\partial}

\newtheorem{proposition}[theorem]{Proposition}
\newcommand {\proof} {\noindent {\bf Proof}. }
\newcommand{\qed}{{ \hfill
                       {\unskip\kern 6pt\penalty 500
                       \raise -2pt\hbox{\vrule\vbox to 6pt{\hrule width 6pt
                       \vfill\hrule}\vrule} \par}   }}

\begin{document}

\title{Stability Analysis of a Simplified Yet Complete Model for Chronic Myelegenous Leukemia}
\author{
Marie Doumic-Jauffret$^{1,2*}$ \and
Peter S.\ Kim$^{3*}$ \and
Beno\^{i}t Perthame$^{1,2}$
}

\maketitle

\begin{abstract}
We analyze the asymptotic behavior of a partial differential equation (PDE) model for hematopoiesis.  This PDE model is derived from the original agent-based model formulated by Roeder {\em et al.} in \cite{roeder2006}, and it describes the progression of blood cell development from the stem cell to the terminally differentiated state.

To conduct our analysis, we start with the PDE model of \cite{kimleelevy2007}, which coincides very well with the simulation results obtained by Roeder \emph{et al}.  We simplify the PDE model to make it amenable to analysis and justify our approximations using numerical simulations.  An analysis of the simplified PDE model proves to exhibit very similar properties to those of the original agent-based model, even if for slightly different parameters.  Hence, the simplified model is of value in understanding the dynamics of hematopoiesis and of chronic myelogenous leukemia, and it presents the advantage of having fewer parameters, which makes comparison with both experimental data and alternative models much easier.
\end{abstract}


\footnotetext[1]{INRIA Paris-Rocquencourt, BANG, BP105, F78153 LeChesnay Cedex. Email: marie.doumic-jauffret@inria.fr  }
\footnotetext[2]{Laboratoire J.-L. Lions, Universit\'{e} Pierre et Marie Curie, CNRS UMR7598, BC187, 4, place Jussieu, F75252 Paris cedex 05. Email: benoit.perthame@upmc.fr}
\footnotetext[3]{Department of Mathematics, University of Utah, Salt Lake City, UT 84112-0090, USA. Email: kim@math.utah.edu}

\renewcommand{\thefootnote}{\fnsymbol{footnote}}
\footnotetext[1]{Contributed comparably to this paper.}

%%%%-----------------------------------------------------------------------------------

\section{Introduction}
\label{section:Introduction}
\setcounter{equation}{0}
\setcounter{figure}{0}
\setcounter{table}{0}

Chronic myelogenous leukemia (CML) is a cancer of the blood and bone marrow that results in the uncontrolled growth of myeloid blood cells.  
More than 90\% of all CML cases are associated with a gene abnormality, known as the Philadelphia (Ph) chromosome \cite{thijsen1999}.
In addition, CML is highly responsive to treatment by the drug imatinib that specifically targets the gene abnormality \cite{druker2000}.

Recently, CML has been the focus of several mathematical models, which we summarize below.
These works were motivated by the desire to explore the mechanisms that control
the disease with the hope that this will lead to new therapeutical strategies.
In 1991, Fokas {\em et al.} formulated the first mathematical model of CML \cite{fokas1991}.
In 2002, Neiman presented a model that accounts for the immune response to CML \cite{neiman2002}.  This work attempted to explain the
transition of leukemia from the stable chronic phase to the accelerated and acute
phases.  A more recent work of Moore and Li in 2004 aimed to identify the
parameters that control cancer remission \cite{moore2004}.  Their main conclusion was that lower growth
rates lead to a greater chance of cancer elimination. In 2005, Komarova {\em et al.} used methods of
stochastic networks to study drug resistance with a particular view toward imatinib \cite{komarova2005}.
In the same year, DeConde et al.
proposed a model for the interaction between the immune system and cancer cells after a
stem cell (or a bone-marrow) transplant \cite{deconde2005}.  The main result of that work was that a slightly
elevated autologous (pretransplant) immune response greatly favors remission.

A new paradigm of cancer development emerged from the idea of cancer stem cells \cite{bonnet1997}.  This hypothesis states that a variety of cancers originate from a self-replenishing, cancer population, now known as cancer stem cells.  
Using this idea, Roeder {\em et al.} developed a mathematical model of CML stem cells \cite{roeder2006}.  In their model, 
leukemia stem cells continually circulate between proliferating and quiescent states.
This formulation contrasts with the alternative paradigm of Michor {\em et al.} in which leukemia cells differentiate progressively from stem cells to differentiated cells without circulating back to previous and more dormant states \cite{michor2005, michor2008}.
Both the Roeder and Michor models are directed to studying the dynamics of imatinib treatment.  However, each model also presents a general paradigm for hematopoiesis that describes blood cell development with or without CML.
In this paper, we focus on the model of Roeder {\em et al.}.

The hypothesis of the Roeder model is based on experimental studies that demonstrated that within a two week period, nearly all hematopoietic stem cells enter cell cycle \cite{bradford1997, cheshier1999}.  A summary of the experimental results and interpretations is also presented in \cite{roeder2005}.
Based on these findings, Roeder {\em et al.} formulated an agent-based model (ABM) to capture probabilistic effects and intrinsic heterogeneity of the stem cell population \cite{roeder2006}.
However, since the ABM is computationally demanding, Roeder {\em et al.} developed an analogous PDE model that is based on their original ABM model \cite{roeder2009}.
In a parallel and independent work, Kim {\em et al.} also developed an analogous PDE model \cite{kimleelevyBMB2008b}.  The main difference between the two PDE models, is that the version of Roeder {\em et al.} simplifies the original formulation by eliminating the explicit representation of the cell cycle, whereas the version of Kim {\em et al.} includes the full complexity of the original ABM.

In this paper, we take the PDE model in \cite{kimleelevyBMB2008b} and simplify it as much as possible without altering the fundamental assumptions of Roeder {\em et al.} \cite{roeder2006}.  
Then we conduct an analysis of the asymptotic behavior of the simplified model for hematopoiesis.
The assumptions of Roeder {\em et al.} include the fact that stem cells exist in two growth environments: proliferating and quiescent.  In addition, proliferating stem cells gradually progress toward further levels of  differentiation until they differentiate completely.  However, proliferating stem cells can reenter the quiescent state, in which they cease dividing and regress toward more primitive levels of differentiation, even up to the fully immature state.  One can interpret that while quiescent stem cells cease dividing, they begin performing other regenerative functions such as returning to the stem cell niche as suggested in \cite{roeder2006}.

The paper is organized as follows.  
In Section~\ref{section:PDEmodel}, we present the PDE model from \cite{kimleelevyBMB2008b}, which proved to coincide very well with the simulation results of Roeder {\em et al.} in \cite{roeder2006}.
In Section~\ref{section:Approximations}, we present four approximations of the PDE model, one of which is very similar to that presented in \cite{roeder2009}.
In Section~\ref{section:NumericalSimulations}, we compare numerical solutions of the four approximations with simulations of the original ABM.
%MD CHANGE 03/15
In Section~\ref{section:Analysis}, we analyze the asymptotic behavior of the simplest models given by Approximations 3 and 4. 
%\ref{section:LinkWithDDE}
%\ref{section:AsymptoticAnalysis}

%%%%-----------------------------------------------------------------------------------

\section{PDE version of  Roeder's model}
\label{section:PDEmodel}
\setcounter{equation}{0}
\setcounter{figure}{0}
\setcounter{table}{0}

\begin{figure}[htbp]
\centering
\includegraphics[scale = 0.45]{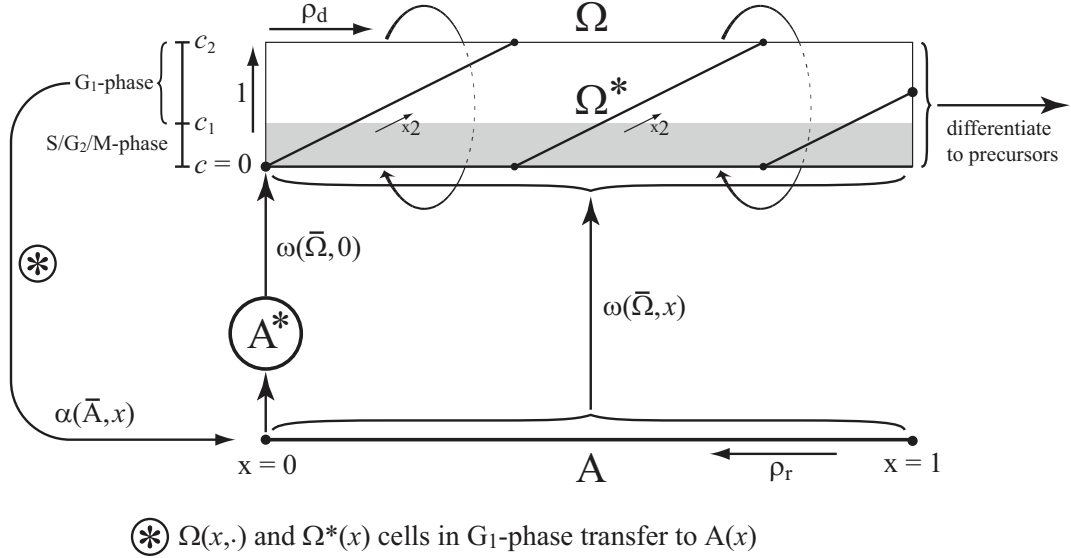}
\caption{State space for the PDE model of CML dynamics in \cite{kimleelevyBMB2008b}.  
The variable $A(x,t)$ represents stem cells in the Alpha (nonproliferating) compartment with intermediate $x$-values, and the variable 
$A^*(t)$ represents cells in the Alpha compartment that have attained the minimum $x$-value of $0$.
The variables $\Omega(x,c,t)$ and $\Omega^*(x,t)$ represent the stem cells in the Omega (proliferating) compartment, where
$\Omega$ and $\Omega^*$ correspond to cells transitioning from the $A$ and $A^*$ subpopulations, respectively.
%The shaded region of $\Omega$ space between $c=0$ and $c=c_1$ corresponds to 
%the S, G$_2$, and M phases of the cell cycle.  The unshaded region between $c=c_1$ 
%and $c=c_2$ corresponds to the G$_1$ phase of the cell cycle.
}
\label{figure:PDEStateSpace}
\end{figure}

In this section, we present a PDE version of the original ABM for hematopoiesis developed by Roeder {\em et al.} in \cite{roeder2006}.
This PDE model is taken from \cite{kimleelevyBMB2008b}.
A state diagram for the model is shown in Figure~\ref{figure:PDEStateSpace}.

In this model, hematopoietic stem cells (HSCs) operate in two growth compartments: 
quiescent (Alpha) and proliferating (Omega).
Stem cells transfer
from Alpha to Omega or from Omega to Alpha at rates $\omega$ and $\alpha$, respectively, where
\begin{align}
\begin{split}
\alpha(x, \Alphatotal) & = e^{-\gamma x} f_{\alpha} (\Alphatotal) , \\
\omega(x, \Omegatotal) & = a_{\text{min}} e^{\gamma x} f_{\omega} (\Omegatotal).
\end{split}
    \label{equation:transitionprobabilities}
\end{align}
The coordinate $x$ is a state variable that ranges from 0 to 1, and
the parameter $a_{\text{min}} = 0.002$ as estimated in \cite{roeder2006}.
The variables $\Alphatotal$ and $\Omegatotal$ denote the population of cells in the Alpha and Omega compartments, respectively, and $f_{\alpha}$ and $f_{\omega}$ 
are sigmoidal functions, whose definition can be  found in Appendix~\ref{appendix:parameterestimates}.

The equations for the PDE model in \cite{kimleelevyBMB2008b} are as follows:
\begin{align}
& \frac{\partial A}{\partial t} - \rho_r \frac{\partial A}{\partial x} 
  = - \omega(x, \Omegatotal) A(x,t)
      + \alpha(x, \Alphatotal) \left( \int_{c_1}^{c_2} \Omega(x,c,t) dc 
      + \mathbf{1}_{\text{G}_1}(x) \Omega^*(x,t) \right)
      \label{equation:PDEalpha} \\
& \frac{dA^*}{dt}
      = \rho_r A(0,t) - \omega(0, \Omegatotal) A^* (t)
      \label{equation:PDEalphastar} \\
& \frac{\partial \Omega}{\partial t} + \rho_d \frac{\partial \Omega}{\partial x} + \frac{\partial \Omega}{\partial c}
  = - \alpha(x, \Alphatotal) \mathbf{1}_{ [c_1, c_2) } (c) \Omega (x,c,t) 
      \label{equation:PDEomega} \\
& \frac{\partial \Omega^*}{\partial t} + \rho_d \frac{\partial \Omega^*}{\partial x}
  = -\alpha(x, \Alphatotal) \mathbf{1}_{\text{G}_1}(x) \Omega^* (x,t)
      \label{equation:PDEomegastar}
\end{align}
with boundary conditions
\begin{align}
& A(1,t) = 0,
    \label{equation:PDEalphaBC} \\
& \Omega (x, 0 ,t) = \Omega (x, c_2 ,t) + \omega(x,\Omegatotal) A(x,t),
    \label{equation:PDEomegaBC1} \\
& \Omega^*(0,t) = \frac{\omega(0,\Omegatotal)}{\rho_d} A^* (t),
    \label{equation:PDEomegastarBC1} \\
& \Omega (x, c_1^+ ,t) = 2 \Omega(x, c_1^- ,t),
    \label{equation:PDEomegaBC2} \\
& \Omega^* ((k \rho_d c_2 + c_1)^+, t) =
    2 \Omega^* ((k \rho_d c_2 + c_1)^- ,t), \qquad k=0,1,2, \ldots .
    \label{equation:PDEomegastarBC2}
\end{align}
In these equations, $x \in [0,1]$, $c \in [0, c_2]$, and
\begin{align*}
\text{G}_1 & = \left( \bigcup_{k=0}^{\infty} \big[ k \rho_d c_2 + c_1, (k+1) \rho_d c_2 \big) \right) \bigcap [0,1].
\end{align*} 
The parameter $c_1$ represents the duration of the combined S/G$_2$/M-phases, and $c_2$ represents the duration of the entire cell cycle. The set $\text{G}_1$ is the set of $x$-values for which $\Omega^*$ cells are in the G$_1$ phase (defined by $c \in [c_1,c_2)$).
In addition, the total Alpha and Omega populations are given by
\begin{align}
\begin{split}
\Alphatotal (t) & = \int_0^1 A(x,t) dx + A^*(t), \\ %\qquad
\Omegatotal (t) & = \int_0^1 
                    \int_0^{c_2} \Omega(x,c,t) dc \, dx
                     + \int_0^1 \Omega^*(x,t) dx.
\end{split}
    \label{equation:PDEtotals}
\end{align}
As shown in the state diagram in Figure~\ref{figure:PDEStateSpace}, the Alpha cells are composed of two subpopulations: $A$ and $A^*$.  
As time progresses, $A$ cells decrease their $x$-coordinate at rate $\rho_r$, until they attain the minimum $x$-value of $0$, at which point they enter $A^*$.

In a similar manner, the Omega cells are divided into two subpopulations: $\Omega$ and $\Omega^*$.
The subpopulation $\Omega$ corresponds to cells that have entered the Omega state from $A$, and the subpopulation $\Omega^*$ corresponds to cells that have entered from $A^*$.
As time progresses, all Omega cells increase their $x$-coordinate at rate $\rho_d$, until they attain the maximum $x$-value of $1$, at which point they differentiate into proliferating precursors.
In addition, all Omega cells progress through the cell cycle as they undergo proliferation.
The G$_1$ phase corresponds to the period during which the cell generates new organelles.
Only cells in the G$_1$ phase can transfer to $A$.
The other phases, S, G$_2$, and M, correspond to the period during which a cell is undergoing mitosis.  Cells in these phases cannot transfer.
A cell's position in the cycle is measured by its $c$-coordinate, which ranges from $0$ to $c_2$.

The first term on the RHS of (\ref{equation:PDEalpha}) accounts for the cells
that transfer out of $A$ into $\Omega$.  The transition rate $\omega$ 
is given by (\ref{equation:transitionprobabilities}),
where $\Omegatotal$ is defined in (\ref{equation:PDEtotals}).  
The second term on the RHS of (\ref{equation:PDEalpha}) is the rate at which cells transfer into $A$ from $\Omega$ and $\Omega^*$. 
The transition rate $\alpha$ is given by (\ref{equation:transitionprobabilities}),and $\Alphatotal$ is given by (\ref{equation:PDEtotals}).
Only $\Omega$ and $\Omega^*$ cells in the $G_1$ phase (i.e., with time counters $c$ between $c_1$ and $c_2$) can transfer into $A$, which explains the boundaries in the integral in (\ref{equation:PDEalpha}) and the indicator function in (\ref{equation:PDEomega}).
In (\ref{equation:PDEalphastar}), the first term on the RHS is the rate at which cells flow from $A$ into $A^*$.  
These $A$ cells flow from the endpoint $x = 0$ into $A^*$.  The second term on the RHS of (\ref{equation:PDEalphastar})
is the rate at which cells flow out of $A^*$ into $\Omega^*$.  Cells coming from $A^*$ enter $\Omega^*$ at the point $(x,c) = (0,0)$.

The expression on the RHS of (\ref{equation:PDEomega}) represents cells in the G$_1$-phase that transfer out of $\Omega$.
The RHS of (\ref{equation:PDEomegastar}) accounts for the rate that cells in the G$_1$-phase flow out of $\Omega^*$ into $A$.
Boundary condition (\ref{equation:PDEalphaBC}) indicates that all Omega cells that attain the maximal value of $x = 1$ commit to differentiation and do not transfer back to $A$.
The first term of boundary condition (\ref{equation:PDEomegaBC1}) indicates that once cells come to the end of a cell cycle at $c = c_2$, they reset their time counters to $c = 0$.
The second term accounts for the rate that cells transfer from $A$ to $\Omega$ at the boundary $c = 0$.  It is the negative of the first term on the RHS of (\ref{equation:PDEalpha}).
Boundary condition (\ref{equation:PDEomegastarBC1}) accounts for the rate that cells transfer from $A^*$ to $\Omega$ at the point $x =0$.  This expression balances the second term on the RHS of (\ref{equation:PDEalphastar}), scaled by the advection rate because it represents a flux. 
Boundary conditions (\ref{equation:PDEomegaBC2}) and (\ref{equation:PDEomegastarBC2}) account for division whenever cells cross the point $c = c_1$ of the cell cycle.
%%%%-----------------------------------------------------------------------------------
%%%%-----------------------------------------------------------------------------------
%%%%-----------------------------------------------------------------------------------
%%%%-----------------------------------------------------------------------------------
%%%%-----------------------------------------------------------------------------------
%%%%-----------------------------------------------------------------------------------
%%%%-----------------------------------------------------------------------------------
\subsection{Rescaled parameters}
%%%%-----------------------------------------------------------------------------------
For convenience, we rescale the parameters from the original ABM \cite{roeder2006} to make them more appropriate for the PDE formulation presented in Section~\ref{section:PDEmodel}.  (The values of the parameters used in \cite{roeder2006} are listed in Table~\ref{table:RoederParameters}.)
First, we measure time in units of days rather than hours.
Then, instead of using the affinity to characterize a stem cell's proclivity to remain in the Alpha (quiescent) state, we replace the affinity variable, $a$, with $x$, which is defined by
$$
a(x) = e^{-\gamma x},\qquad \gamma = -\log a_{\text{min}} = 6.2146.
$$
Since the affinity ranges from $a_{\text{min}} = 0.002$ to $a_{\text{max}} = 1$, $x$ ranges from 0 to 1.
\\
In \cite{roeder2006}, at every time step of 1 hour, the affinity increases by a factor $r$ for cells in Alpha and decreases by a factor of $d$ for cells in Omega.
Hence, the rescaled advection rates are, taking into account the time unit of days,
\begin{align*}
\rho_r = \frac{24 \log r}{\gamma} = 0.3681, \qquad \rho_d = \frac{24 \log d}{\gamma} = 0.1884.
\end{align*}
We also rescale the cell populations down by a factor of $10^5$ cells, so that the scaling factors $\tilde{N}_A$ and $\tilde{N}_{\Omega}$ become equal to $1$ rather than to $10^5$.
Unlike \cite{roeder2006}, we write the transition rates, $\alpha(\Alphatotal(t), x)$ and $\omega(\Omegatotal(t), x)$ as functions of $x$ and not of $a$. Finally, 
since time is measured in days, we need to multiply the original transition rates $f_{\alpha/\omega}$ by 24 as shown in equation (\ref{equation:sigmoidfunctions}).
Our rescaled parameters are listed in Table~\ref{table:RoederParameters} alongside the original values from \cite{roeder2006}.
%%%%-----------------------------------------------------------------------------------
%%%%-----------------------------------------------------------------------------------
%%%%-----------------------------------------------------------------------------------
%%%%-----------------------------------------------------------------------------------
%%%%-----------------------------------------------------------------------------------
%%%%-----------------------------------------------------------------------------------
%%%%-----------------------------------------------------------------------------------
\section{Approximations of the PDE system}
\label{section:Approximations}
\setcounter{equation}{0}
\setcounter{figure}{0}
\setcounter{table}{0}
%%%%-----------------------------------------------------------------------------------
%%%%-----------------------------------------------------------------------------------
The PDE system of Section~\ref{section:PDEmodel} is rather complex and hence difficult to analyze mathematically.  Thus, we introduce several approximation steps that can be used to simplify the system.
The approximation steps are as follows:
\begin{itemize}
\item[{\bf (0)}] No approximation.  This label refers to the full PDE model 
\eqref{equation:PDEalpha}--\eqref{equation:PDEomegastarBC2} presented in Section~\ref{section:PDEmodel}.
\item[{\bf (1)}] Eliminate the explicit representation of the cell cycle by averaging the system over the variable $c$.
\item[{\bf (2)}] (a) Assume the transition rate $\omega(x,\Omegatotal)$ vanishes for all $x$, except $x = 0$.  In other words, Alpha cells do not transition into Omega, except at the point $(x,c) = (0,0)$.
This approximation causes the variable $\Omega$ vanish. \\
(b) Assume the advection rate $\rho_r$ is effectively infinite so that cells transfering from Omega into Alpha immediately enter $A^*$ without passing through $A$.  This approximation causes the variable $A$ to vanish.
\item[{\bf (3)}] Combine Approximations 1 and 2.
\item[{\bf (4)}] Combine Approximations 1 and 2 with the further approximation that the transition rates $\alpha(x,\Alphatotal)$ and $\Omegatotal(x,\Omegatotal)$ do not depend on $x$.  %In other words, these values are constant for all $x$.
\end{itemize}
%%%%-----------------------------------------------------------------------------------
%%%%-----------------------------------------------------------------------------------
%%%%-----------------------------------------------------------------------------------
%%%%-----------------------------------------------------------------------------------
\subsection{PDE system for Approximation 1}
\label{section:Approximation1}
%%%%-----------------------------------------------------------------------------------
In this approximation, we eliminate the explicit representation of the cell cycle by averaging the system over the variable $c$.
which is  removed, allowing us to combine $\Omega$ and $\Omega^*$ into one single population.
The new unknown $\Omega$ now represents the total Omega population, it satisfies a new PDE that 
replaces equations \eqref{equation:PDEomega} and \eqref{equation:PDEomegastar}, namely
\begin{align}
\frac{\partial \Omega}{\partial t} + \rho_d \frac{\partial \Omega}{\partial x}
& = \omega(x, \Omegatotal) A(x,t) + (-\kappa \alpha(x, \Alphatotal) + b) \Omega(x,t).
    \label{equation:Approx1omega}
\end{align}
This averaging procedure leads to continuous rates of division and transfer, $b$ and $\kappa \alpha(\Alphatotal, x)$, where $b$ is the average growth rate
and $\kappa$ is the proportion of time a cell spends in the G$_1$-phase of the cell cycle.
These parameters are estimated in Section~\ref{section:ApproximationParameters}.
\\
The boundary condition for $\Omega$ at $x = 0$ is the same as the original boundary condition for $\Omega^*$ given by (\ref{equation:PDEomegastarBC1}),
i.e.,
\begin{align}\label{equation:Approx1omegaBC}
\Omega(0,t) & = \frac{\omega(0, \Omegatotal)}{\rho_d} A^* (t), \qquad \text{where }\Omegatotal(t)= \int_0^1 \Omega(x,t)dx.
\end{align}
%%%%-----------------------------------------------------------------------------------

We complete Equations \eqref{equation:Approx1omega} and \eqref{equation:Approx1omegaBC} with the following equations for $A$ and $A^*.$
\begin{align}
& \frac{\partial A}{\partial t} - \rho_r \frac{\partial A}{\partial x} 
  = - \omega(x, \Omegatotal) A(x,t)
      + \kappa \alpha(x, \Alphatotal) \Omega(x,t),  \label{equation:Approx1alpha}
\\
& \frac{dA^*}{dt}
      = \rho_r A(0,t) - \omega(0, \Omegatotal) A^* (t), \label{equation:Approx1alphastar}
\\
& A(1,t) = 0. \label{equation:Approx1alphaBC}
\end{align}
To obtain the PDE \eqref{equation:Approx1alpha}  for $A$ we have just modified (\ref{equation:PDEalpha}) to take into account the changes in the model for the Omega cells. The last term on the RHS balances the corresponding term in (\ref{equation:Approx1omega}).
The ODE for $A^*$ and boundary condition for $A$ remain unchanged (equations \eqref{equation:Approx1alphastar} and \eqref{equation:Approx1alphaBC} are exactly equations (\ref{equation:PDEalphastar}) and (\ref{equation:PDEalphaBC}). 
Approximation 1 is essentially the same PDE system used by Roeder {\em et al.} in \cite{roeder2009}.
%%%%-----------------------------------------------------------------------------------
%%%%-----------------------------------------------------------------------------------
%%%%-----------------------------------------------------------------------------------
%%%%-----------------------------------------------------------------------------------
\subsection{PDE system for Approximation 2}
\label{section:Approximation2}
%%%%-----------------------------------------------------------------------------------
%%%%-----------------------------------------------------------------------------------
In Approximation 2, we suppose $A$ and $\Omega$ are negligible compared to $A^*$ and $\Omega^*$ .  Indeed, numerical simulations from \cite{kimleelevyBMB2008b} show that over 98\% of Alpha cells remain in the $A^*$ compartment and over 94\% of Omega cells remain in $\Omega^*$ compartment over time.  Hence, to a good approximation, we can exclude the $A$ and $\Omega$ populations.
As a result, we are only left with the following equations
\begin{align}
& \frac{dA^*}{dt} = -\omega(0, \Omegatotal) A^*(t) + \int_0^1 \alpha(x, A^*) \mathbf{1}_{\text{G}_1}(x) \Omega^*(x,t) dx , \label{equation:Approx2astar}\\
& \frac{\partial \Omega^*}{\partial t} + \rho_d \frac{\partial \Omega^*}{\partial x}
= -\alpha(x, A^*) \mathbf{1}_{\text{G}_1}(x) \Omega^*(x,t), \label{equation:Approx2omegastar}
\end{align}
with boundary conditions
\begin{align}
& \Omega^*(0,t) = \frac{\omega(0, \Omegatotal)}{\rho_d} A^* (t),  \Omegatotal(t)= \int_0^1 \Omega^*(x,t)dx, \label{equation:Approx2omegastarBC1}\\
& \Omega^* (\rho_d (k c_2 + c_1)^+, t) = 2 \Omega^* ( \rho_d (k c_2 + c_1)^- ,t), \qquad k=0,1,2, \ldots . \label{equation:Approx2omegastarBC2}
\end{align}
The first term on the RHS of \eqref{equation:Approx2astar} for $A^*$ accounts for the rate that cells flow from $A^*$ into $\Omega^*$, and the second term on the RHS accounts for the rate that cells transfer from $\Omega^*$ directly into $A^*$ without passing through $A$. 
In the equations above, \eqref{equation:Approx2astar} replaces \eqref{equation:PDEalpha}, \eqref{equation:PDEalphastar} and \eqref{equation:PDEalphaBC} and can be obtained by integrating \eqref{equation:PDEalpha} with respect to $x$, adding it to \eqref{equation:PDEalphastar}, and assuming that
$$A^* \gg \int_0^1 A(t,x) dx
\hspace{0.5cm}
\text{and}
\hspace{0.5cm}
\int_0^1 \alpha(x,\Alphatotal) \mathbf{1}_{\text{G}_1}(x) \Omega^*(t,x) dx \gg \int_0^1 \int_{c_1}^{c_2} \alpha(x,\Alphatotal) \Omega (t,x,c) dx dc.$$
Equation
\eqref{equation:Approx2omegastar}
and boundary conditions 
\eqref{equation:Approx2omegastarBC1} and \eqref{equation:Approx2omegastarBC2}
for $\Omega^*$ remain the same as 
(\ref{equation:PDEomegastar}), (\ref{equation:PDEomegastarBC1}) and
(\ref{equation:PDEomegastarBC2}), respectively.
%%%%-----------------------------------------------------------------------------------
%%%%---------------------------------------------------------------------------%%%%-----------------------------------------------------------------------------------
%%%%-----------------------------------------------------------------------------------
%%%%-----------------------------------------------------------------------------------
\subsection{PDE system for Approximation 3}
%\label{section:Approximation3}
%%%%-----------------------------------------------------------------------------------

We combine Approximations 1 and 2 
and to simplify notation, replace $\kappa \alpha(x,\Alphatotal)$ by $\alpha(x,A^*)$ and $ \omega(0, \Omegatotal)$ by $\omega(\Omegatotal)$.
Then we arrive at the system
\begin{align}
& \frac{dA^*}{dt} = -\omega(\Omegatotal) A^*(t) + %\kappa 
\int_0^1 \alpha(x, A^*) \Omega^*(x,t) dx , \label{equation:Approx3alphastar}
\\
& \frac{\partial \Omega^*}{\partial t} + \rho_d \frac{\partial \Omega^*}{\partial x}
  = \left( -
%\kappa 
\alpha(x, A^*) + b \right) \Omega^*(x,t) \label{equation:Approx3omegastar}
\end{align}
with boundary condition \eqref{equation:PDEomegastarBC1}, i.e.,
\be
\label{equation:Approx3omegastarBC}
\Omega^*(0,t) = \frac{\omega(\Omegatotal)}{\rho_d} A^* (t), \qquad \Omegatotal(t) = \int_0^1 \Omega^*(x,t)dx.
\ee
The parameter $b$ 
%and $\kappa$ are 
is defined in the same way as in Section~\ref{section:Approximation1} and is estimated in Section~\ref{section:ApproximationParameters}.
%%%%-----------------------------------------------------------------------------------
%%%%-----------------------------------------------------------------------------------
%%%%-----------------------------------------------------------------------------------
%-------------------------------------------------------------------
\subsection{PDE system for Approximation 4}
\label{section:Approximation4}
%%%%-----------------------------------------------------------------------------------
We combine Approximations 1 and 2 with the further approximation that the transition rate $\alpha(x,\Alphatotal)$  does not depend on $x$. To do so we may replace $\alpha(x, \Alphatotal)$ with its midpoint approximation 
$ \alpha(A^*) =\alpha(0.5, A^*)$.  In addition, we  write $\omega(\Omegatotal)$ in place of $\omega(0, \Omegatotal)$.
Thus, we arrive at the simplified system
\begin{align}
& \frac{dA^*}{dt} = -\omega(\Omegatotal) A^*(t) + \alpha(A^*) \int_0^1 \Omega^*(x,t) dx , 
\label{equation:Approx4alphastar} \\
& \frac{\partial \Omega^*}{\partial t} + \rho_d \frac{\partial \Omega^*}{\partial x}
  = \left( -\alpha(A^*) + b \right) \Omega^*(x,t), \label{equation:Approx4omegastar}
\end{align}
with boundary condition
\begin{align}
\Omega^*(0,t) & = \frac{\omega(\Omegatotal)}{\rho_d} A^* (t),  \qquad \Omegatotal(t) = \int_0^1 \Omega^*(x,t)dx.
\label{equation:Approx4omegastarBC}
\end{align}
%MD CHANGE 15/03: I erase the following sentence which is already written at the end of the previous section: 
%The parameters $b$ and $\kappa$ are defined in the same way as in Section~\ref{section:Approximation1} and are estimated in Section~\ref{section:ApproximationParameters}.
%%%%-----------------------------------------------------------------------------------
%%%%-----------------------------------------------------------------------------------
\subsection{Additional parameters}
\label{section:ApproximationParameters}
%%%%-----------------------------------------------------------------------------------
%%%%-----------------------------------------------------------------------------------
The parameters $b$ and $\kappa$ depend on how frequently cells circulate back and forth between the Alpha and Omega compartments.  In particular, frequent turnover results in a decreased cell cycle.
In \cite{roeder2009}, Roeder {\em et al.} estimate $b$ and $\kappa$ as functions of the total cell population, $\Alphatotal + \Omegatotal$.  For simplicity, we assume the total cell population remains constant at the value $1 \times 10^5$.  Then using the functions in \cite{roeder2009}, it follows that Omega cells divide once every 39.5 hours on average and that they spend 54\% of their time in the G$_1$-phase \cite{roeder2009}.
Hence, we estimate
$$
b = \frac{\log(2)}{c_1+c_2}=\frac{\log(2)}{39.5\text{ hours}} = 0.42\text{/days}
\hspace{1cm}
\text{and}
\hspace{1cm}
\kappa = \frac{c_1}{c_1+c_2}=0.54
$$
These estimates are listed in Table~\ref{table:RoederParameters}.
%%%%-----------------------------------------------------------------------------------
%--------------------------------------------------------------------------------------
\section{Numerical simulations}
\label{section:NumericalSimulations}
\setcounter{equation}{0}
\setcounter{figure}{0}
\setcounter{table}{0}
%%%%-----------------------------------------------------------------------------------
%%%%-----------------------------------------------------------------------------------
We compare numerical solutions of the four PDE models presented in Section~\ref{section:Approximations} to simulations of the ABM from \cite{roeder2009}.  
For the transport PDEs we use the explicit upwind scheme while the ABM uses the algorithm from \cite{roeder2006}, which is summarized in Appendix~\ref{appendix:RoederAlgorithm}.
All programs are run in Matlab 7.6.0.

\begin{figure}[htbp]
\begin{center}
\begin{tabular}{ccc}
\includegraphics[scale = 0.39]{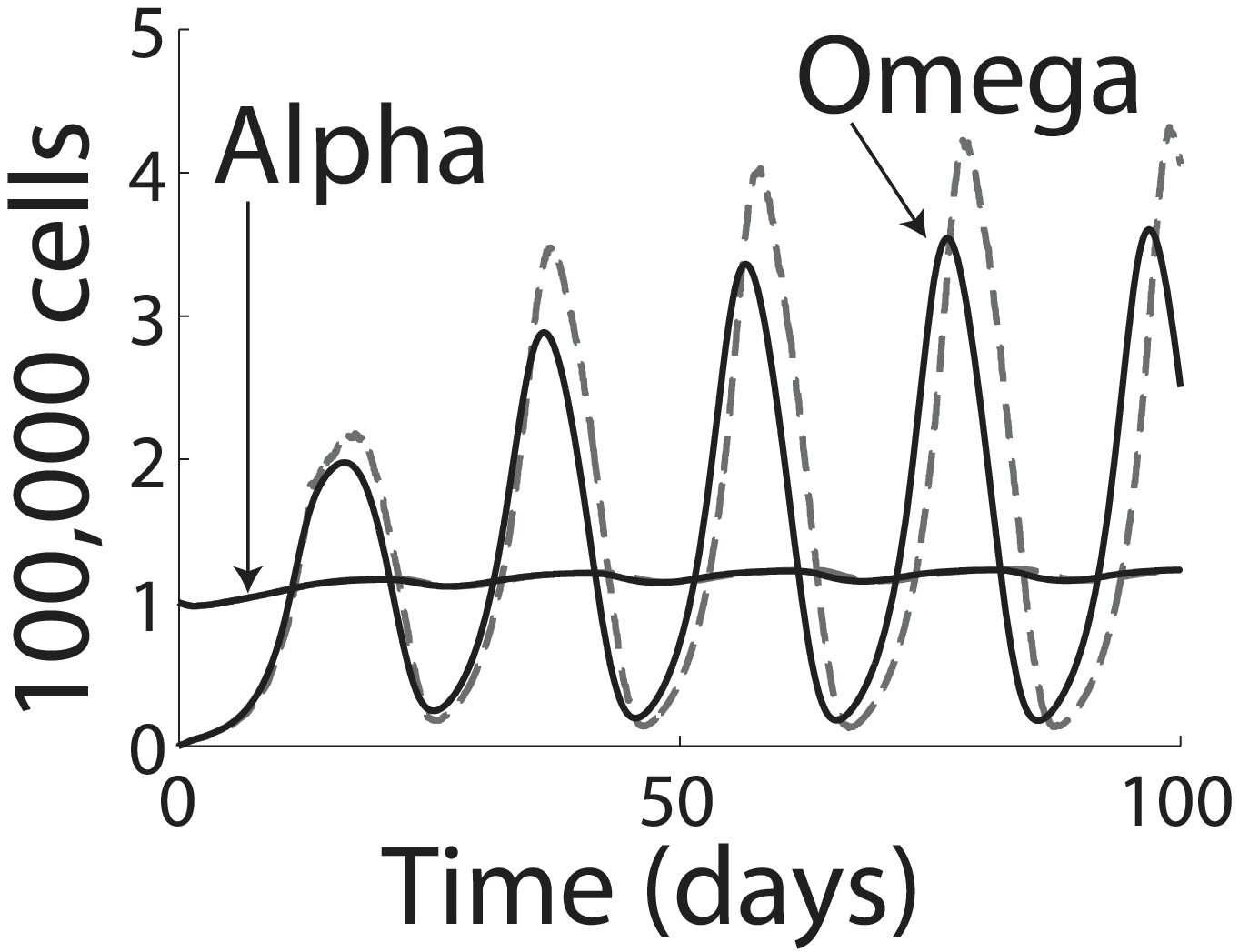} &
    \includegraphics[scale = 0.39]{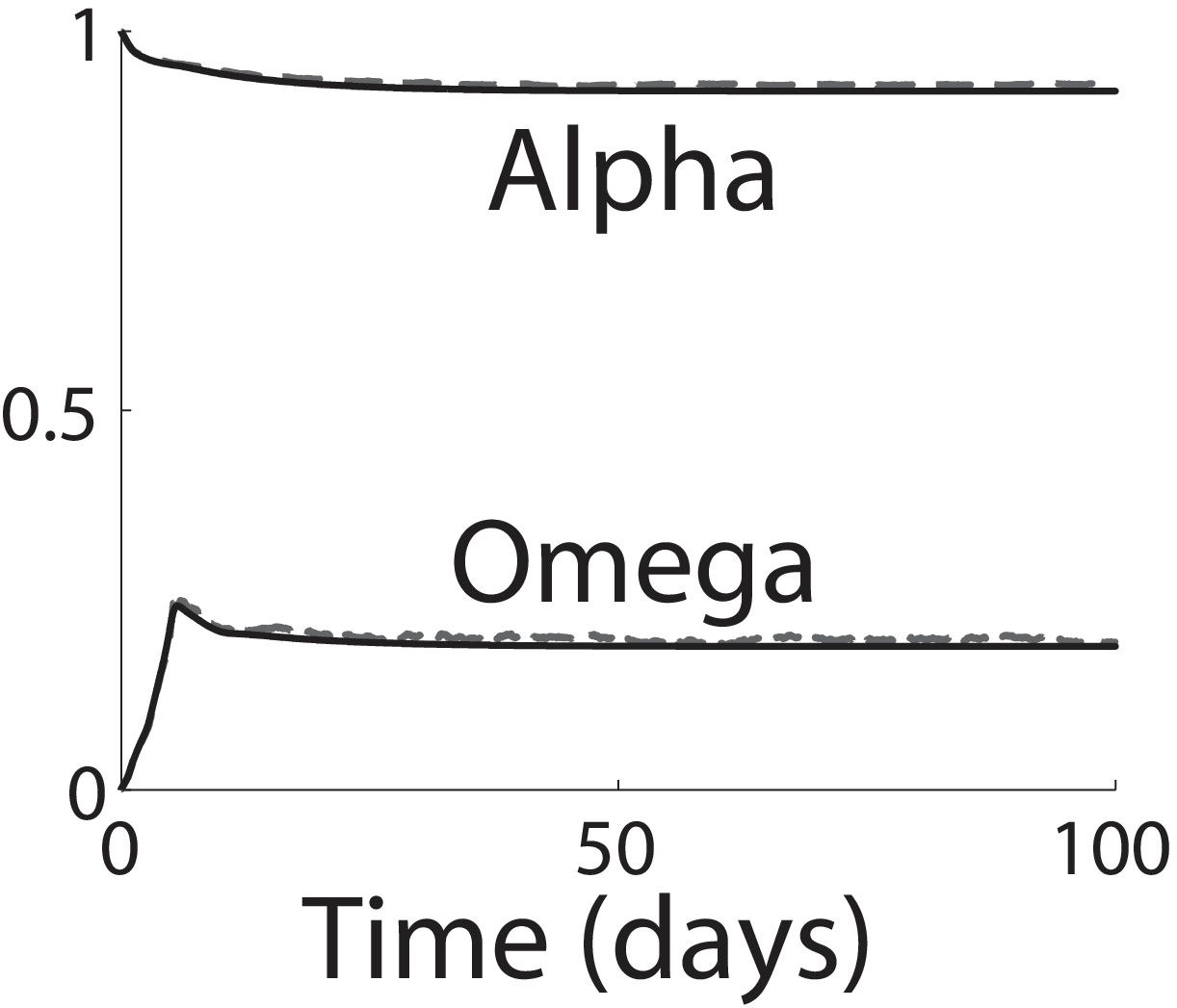} &
    \includegraphics[scale = 0.39]{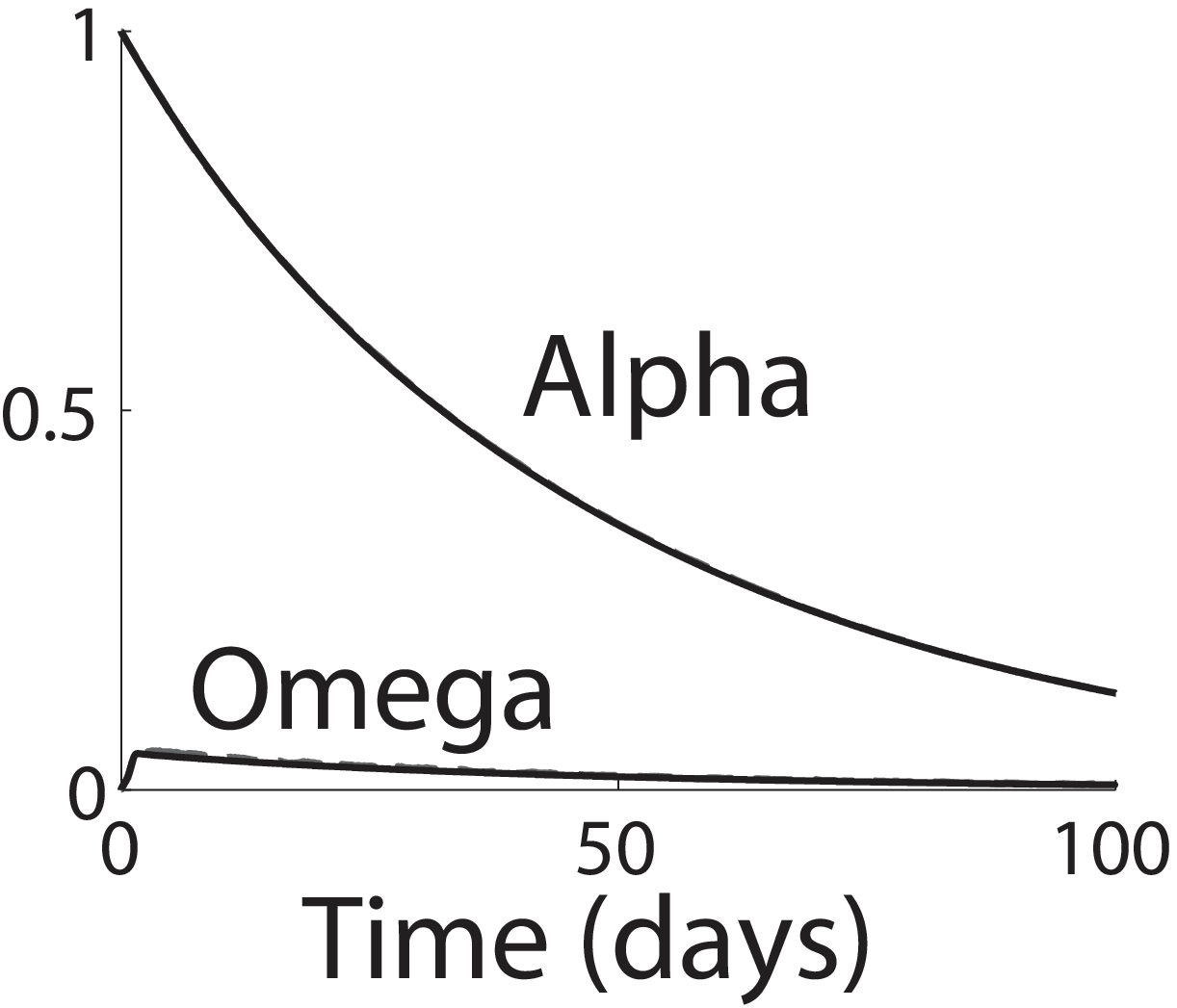} \\
{\bf (a)} & {\bf (b)} & {\bf (c)} \\
\end{tabular}
\end{center}
\caption{
Time plots of numerical solutions to Approximation 0 and simulations of the ABM.
Solutions to Approximation 0 are shown by solid black lines, and
simulations of the ABM are shown by dashed gray lines.  The results of Approximation 0 and the ABM are almost indistinguishable in all cases, except for the Omega population in plot (a).
(a) Time plots for $d = 1.02$.
(b) Time plots for $d = 1.05$.
(c) Time plots for $d = 1.2$.
}
\label{figure:ComparisonApprox0toABM}
\end{figure}

Figure~\ref{figure:ComparisonApprox0toABM} shows time plots of numerical solutions to Approximation 0 (i.e., the full PDE model) and simulations of the ABM for 100 days.  In each of the plots, the initial condition is $A^*(t=0) = 1$ while all other variables are 0.
All parameters are taken from Table~\ref{table:RoederParameters}, except in Figures~\ref{figure:ComparisonApprox0toABM}a, b, and c, where the differentiation rate $d = 1.02$, 1.05, and 1.2, respectively.  (For the PDE model, these values of $d$ yield $\rho_d$ values of 0.0765, 0.1884, and 0.7041, respectively.)
%%%%-----------------------------------------------------------------------------------

The plots in Figures~\ref{figure:ComparisonApprox0toABM}a, b, and c show the following characteristic behaviors, respectively: approach to periodic behavior,
approach to a nonzero steady state, and approach to the zero steady state.
Furthermore, we see that the behavior of the full PDE model closely matches the behavior of the ABM in all three cases.

\begin{figure}[htbp]
\begin{center}
\begin{tabular}{ccc}
\includegraphics[scale = 0.39]{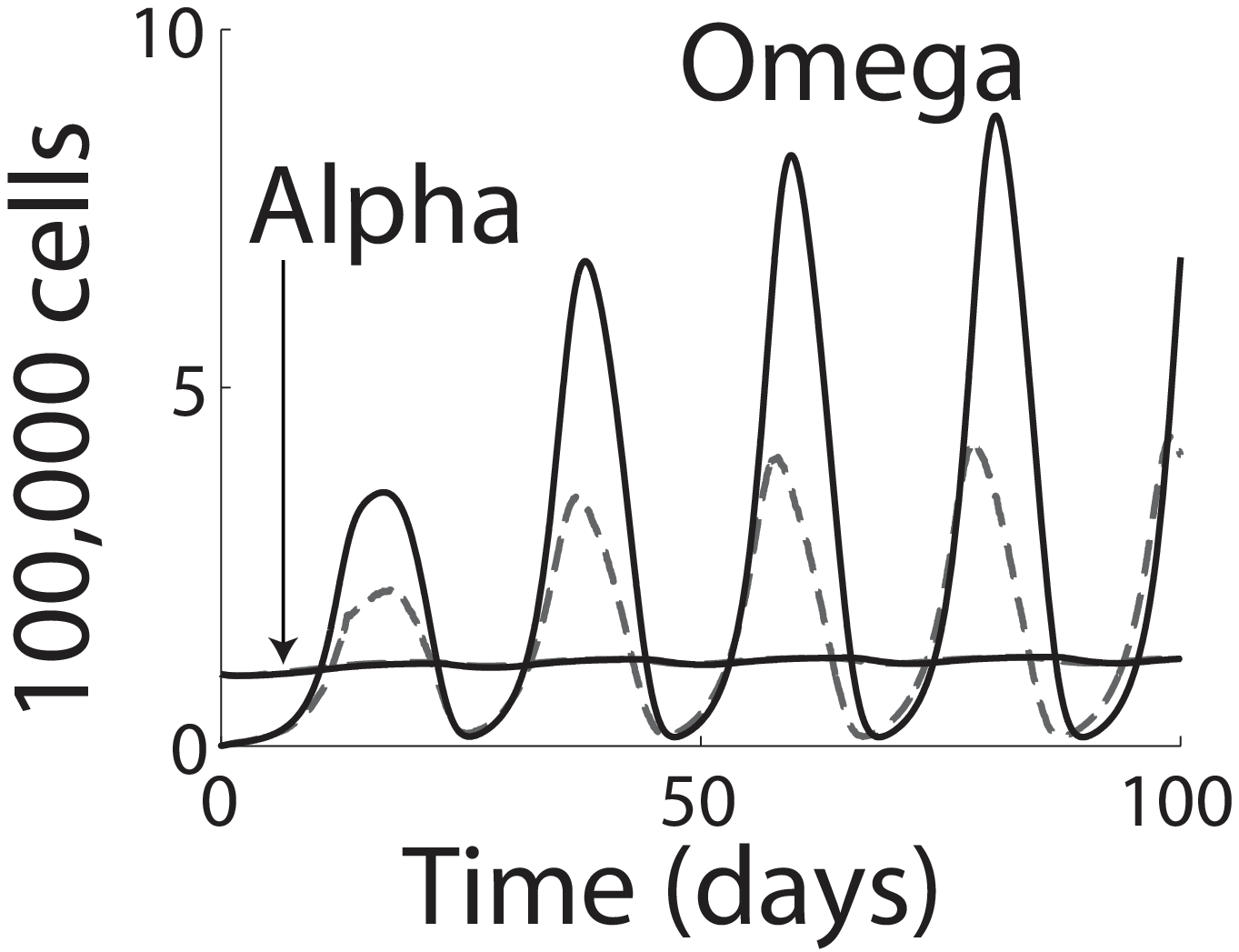} &
    \includegraphics[scale = 0.39]{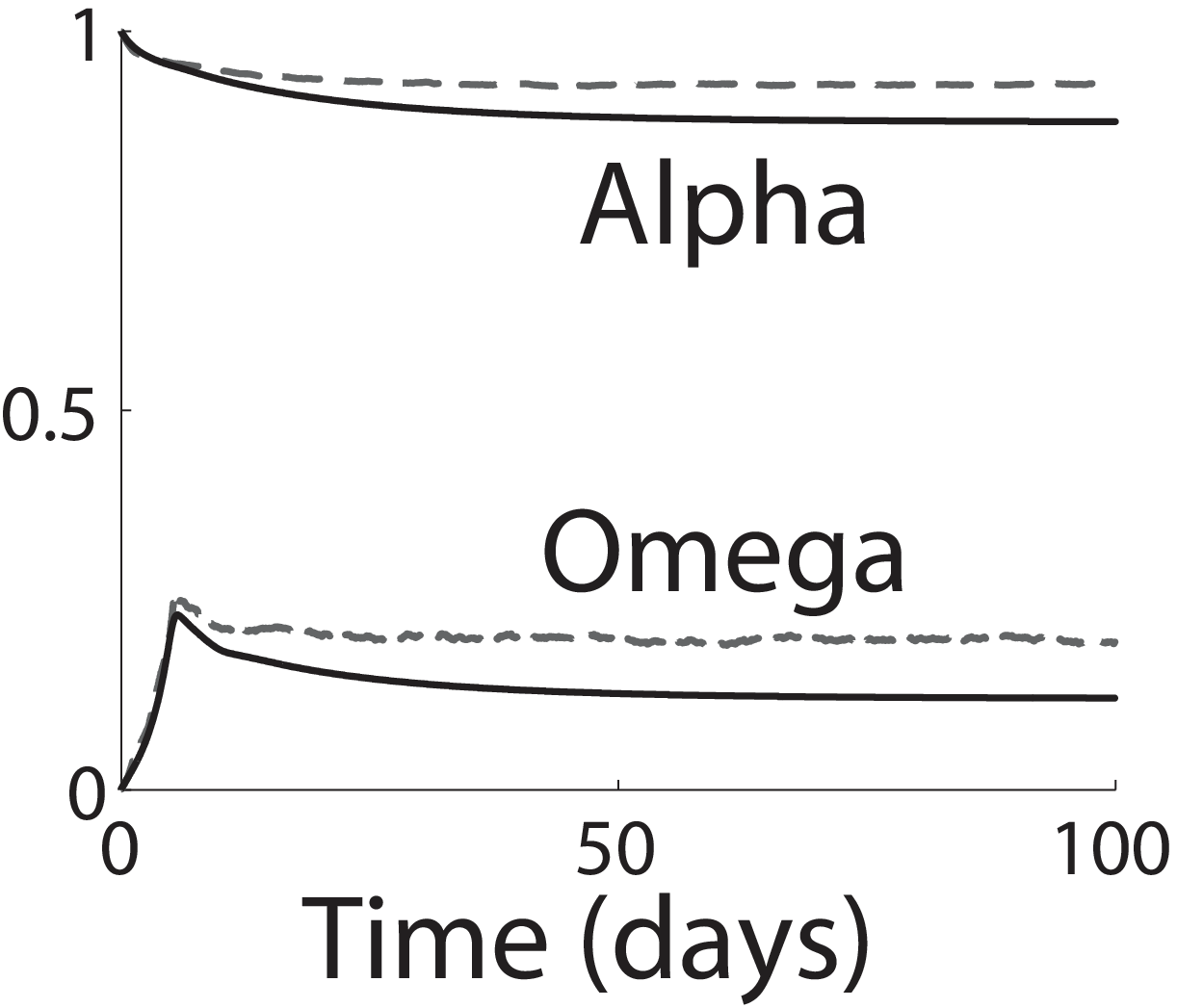} &
    \includegraphics[scale = 0.39]{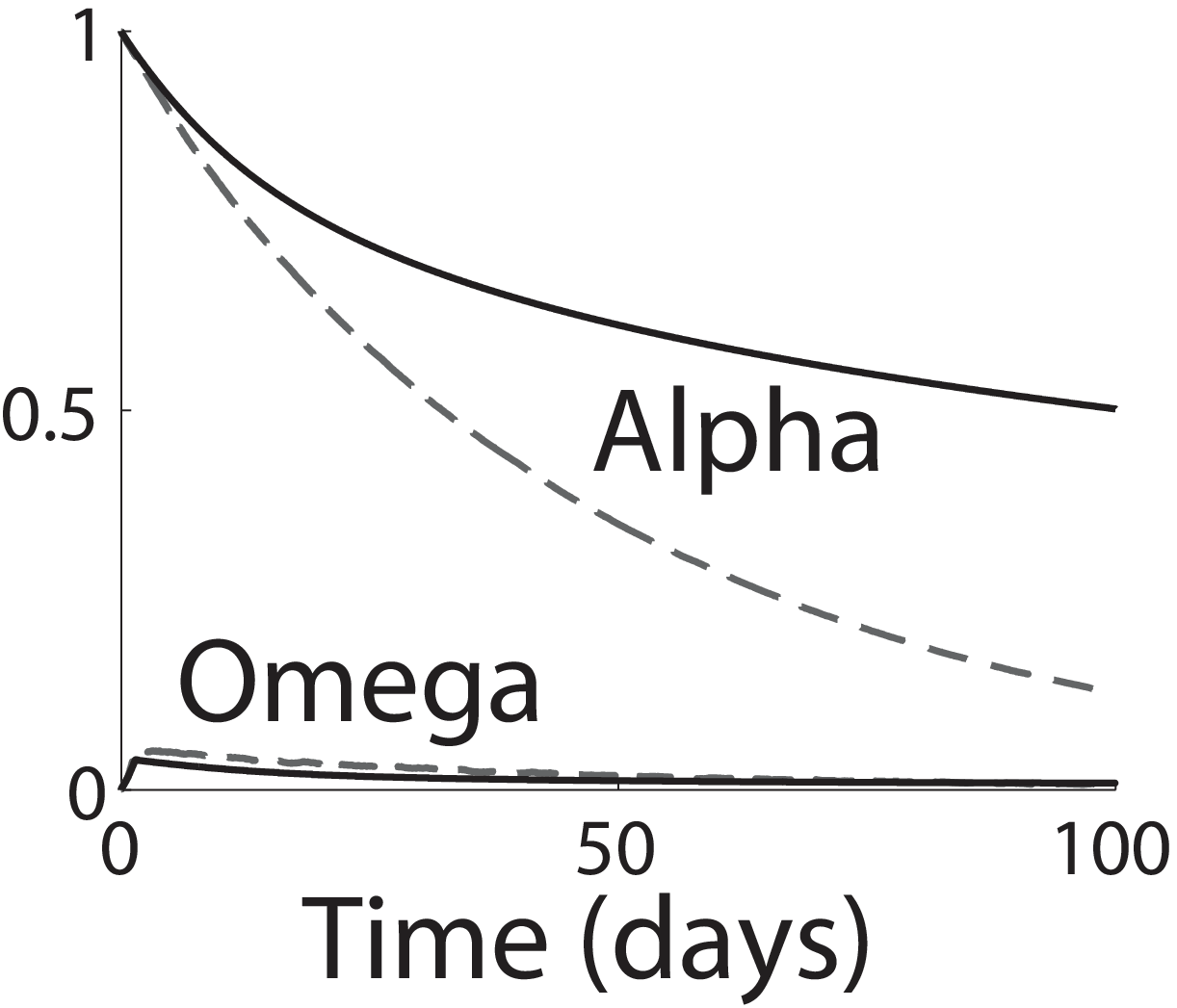} \\
{\bf (a)} & {\bf (b)} & {\bf (c)} \\
\end{tabular}
\end{center}
\caption{
Time plots of numerical solutions to Approximation 1 and simulations of the ABM.
Solutions to Approximation 1 are shown by solid black lines, and
simulations of the ABM are shown by dashed gray lines.
(a) Time plots for $d = 1.02$.
(b) Time plots for $d = 1.05$.
(c) Time plots for $d = 1.2$.
}
\label{figure:ComparisonApprox1toABM}
\end{figure}

Using the same initial conditions and parameters as above,
Figure~\ref{figure:ComparisonApprox1toABM} compares Approximation 1 to the ABM.
We see that the behavior of Approximation 1 qualitatively matches the behavior of the ABM to a large extent in all three cases.
%%%%-----------------------------------------------------------------------------------

\begin{figure}[htbp]
\begin{center}
\begin{tabular}{ccc}
\includegraphics[scale = 0.39]{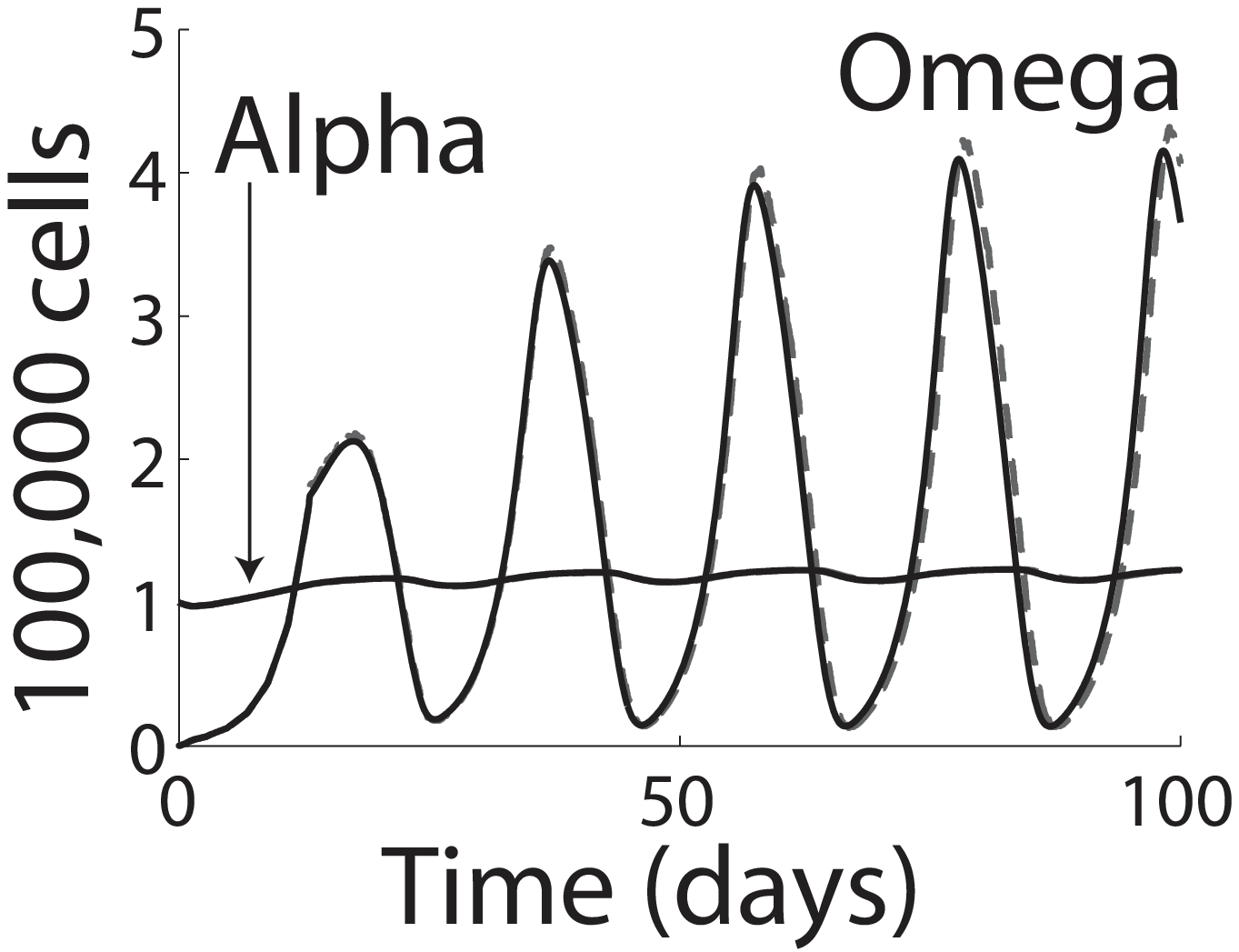} &
    \includegraphics[scale = 0.39]{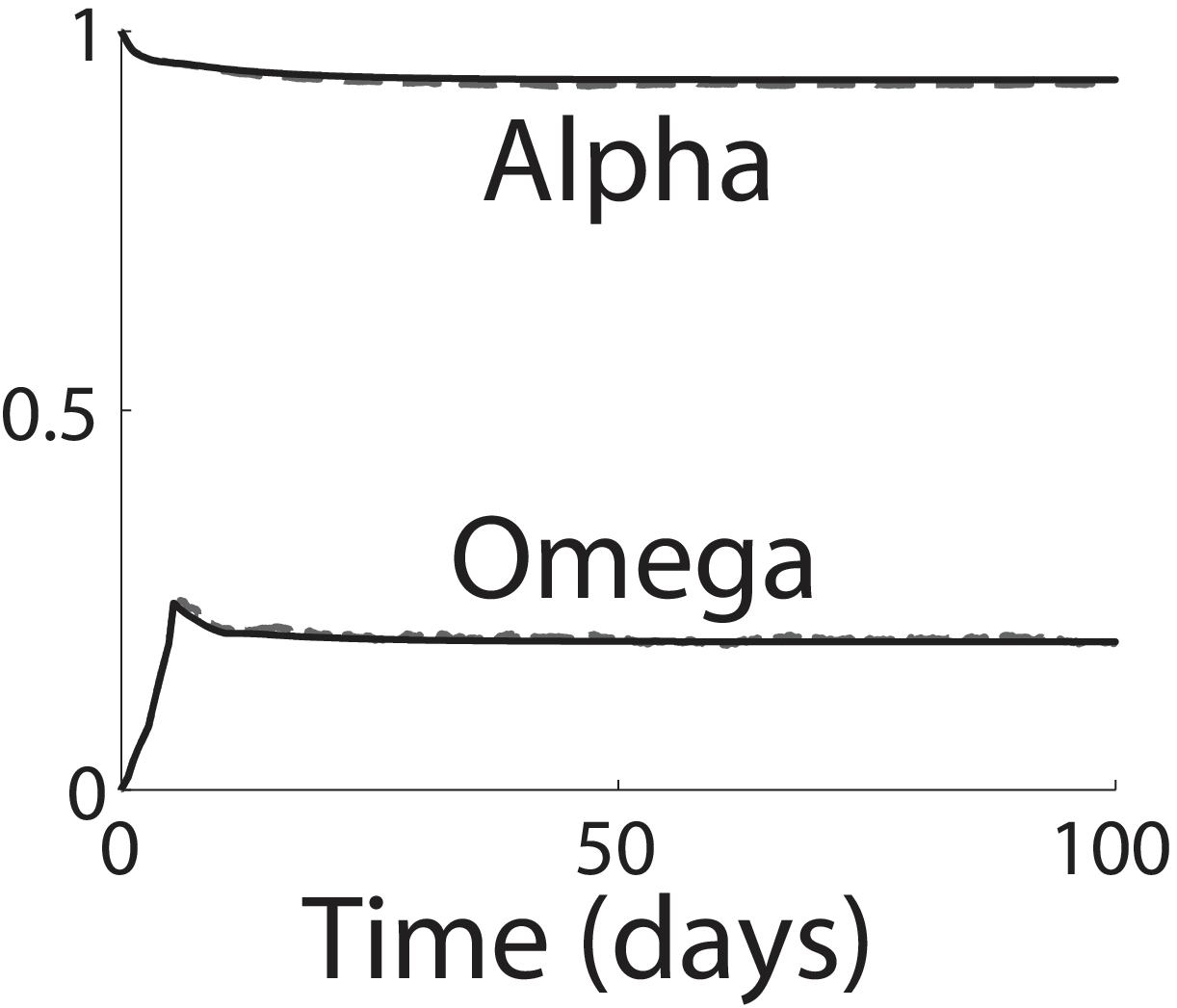} &
    \includegraphics[scale = 0.39]{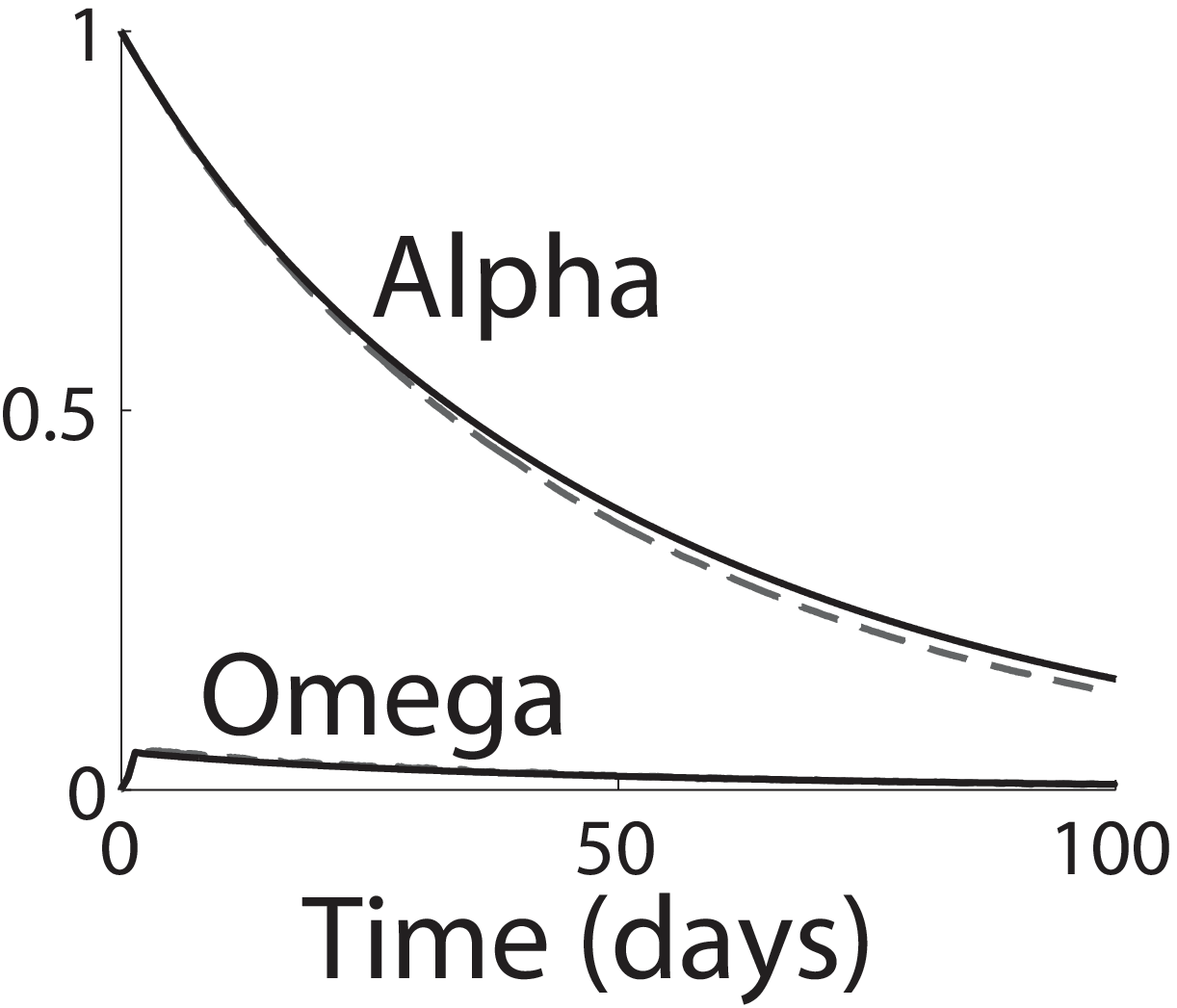} \\
{\bf (a)} & {\bf (b)} & {\bf (c)} \\
\end{tabular}
\end{center}
\caption{
Time plots of numerical solutions to Approximation 2 and simulations of the ABM.
Solutions to Approximation 2 are shown by solid black lines, and
simulations of the ABM are shown by dashed gray lines.
(a) Time plots for $d = 1.02$.
(b) Time plots for $d = 1.05$.
(c) Time plots for $d = 1.2$.
}
\label{figure:ComparisonApprox2toABM}
\end{figure}

Using the same initial conditions and parameters as before,
Figure~\ref{figure:ComparisonApprox2toABM} compares Approximation 2 and the ABM.
We see that the behavior of Approximation 2 closely replicates the behavior of the ABM in all three cases.  In fact, the results of the two models are almost indistinguishable in all cases.
%%%%-----------------------------------------------------------------------------------

\begin{figure}[htbp]
\begin{center}
\begin{tabular}{ccc}
\includegraphics[scale = 0.39]{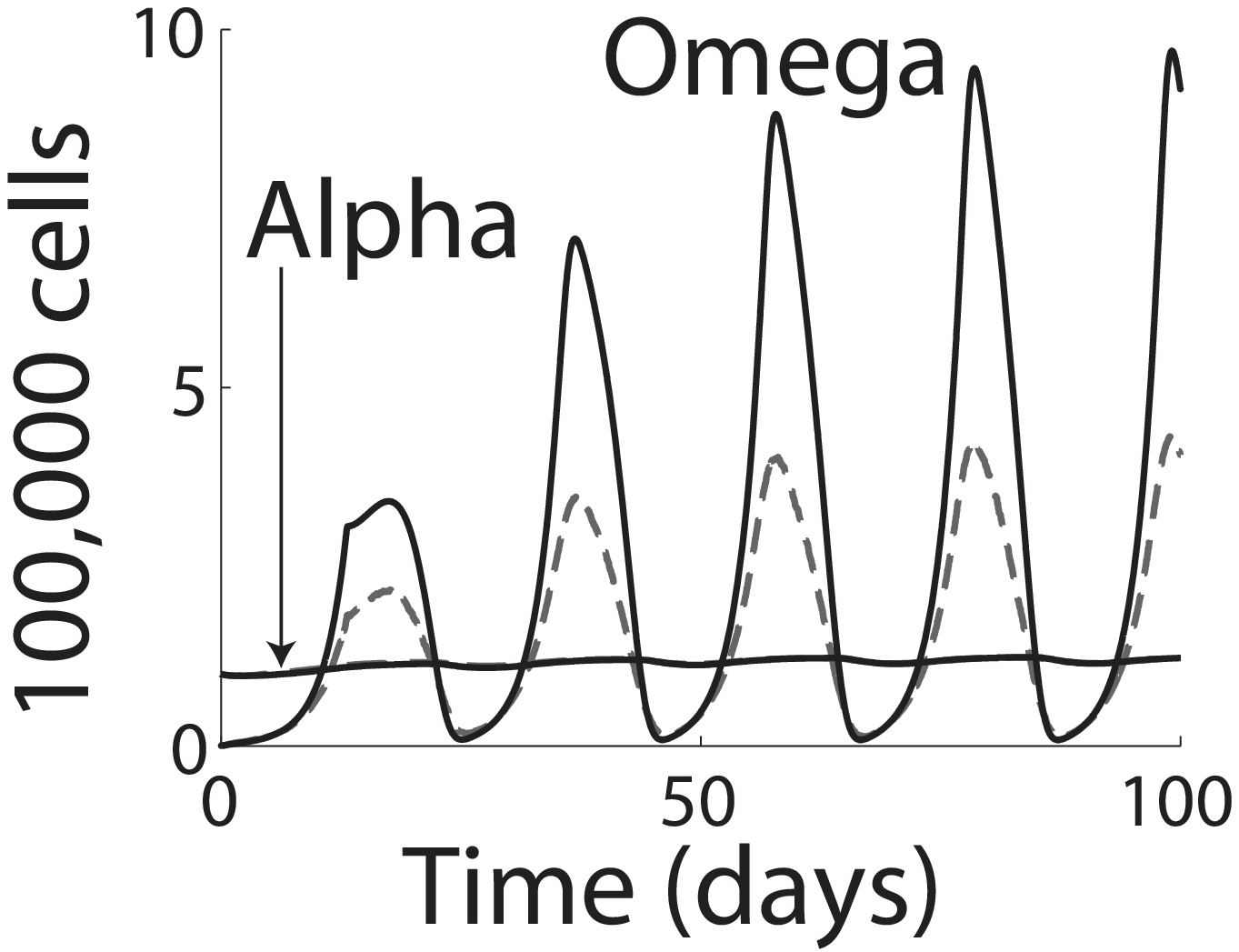} &
    \includegraphics[scale = 0.39]{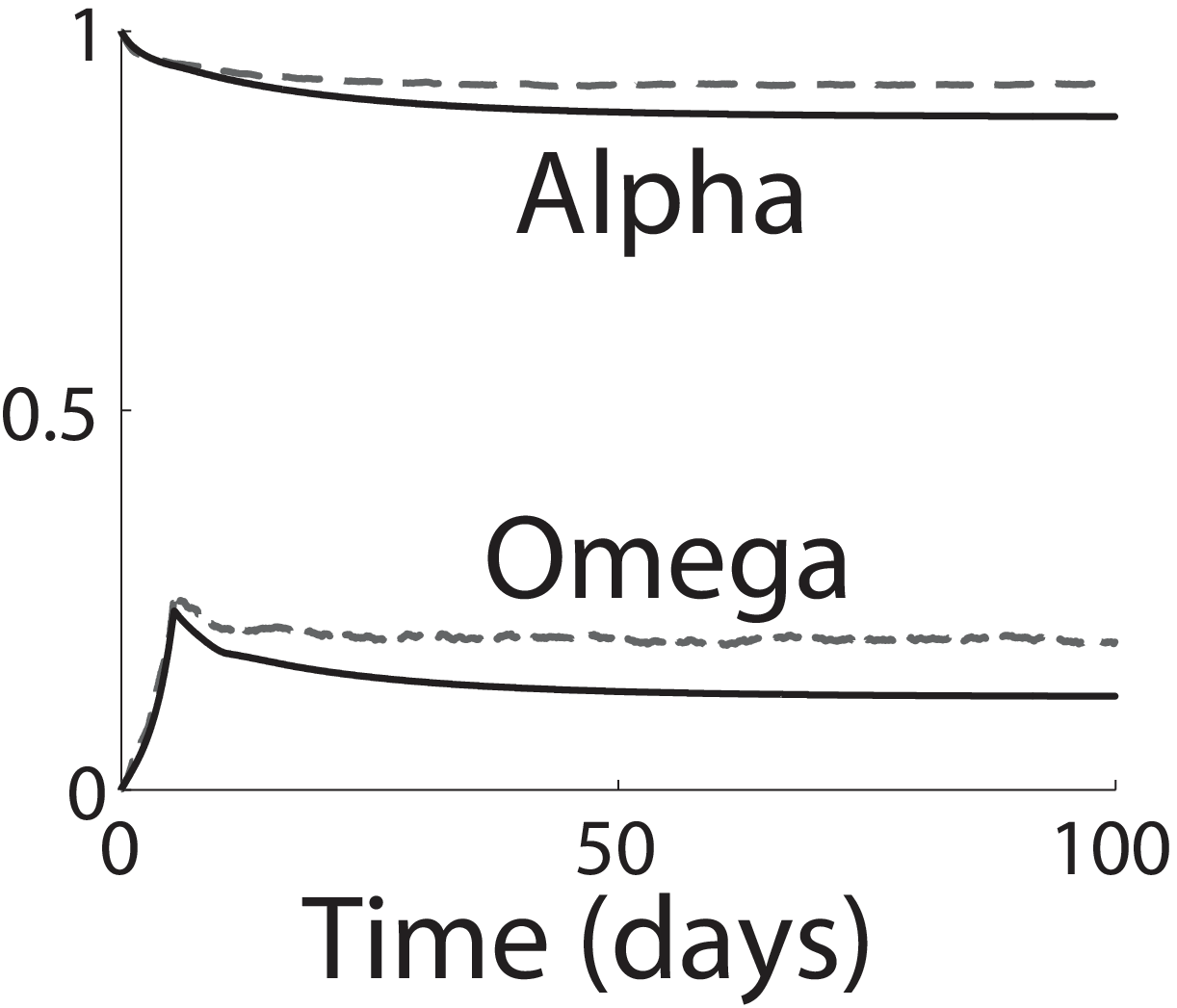} &
    \includegraphics[scale = 0.39]{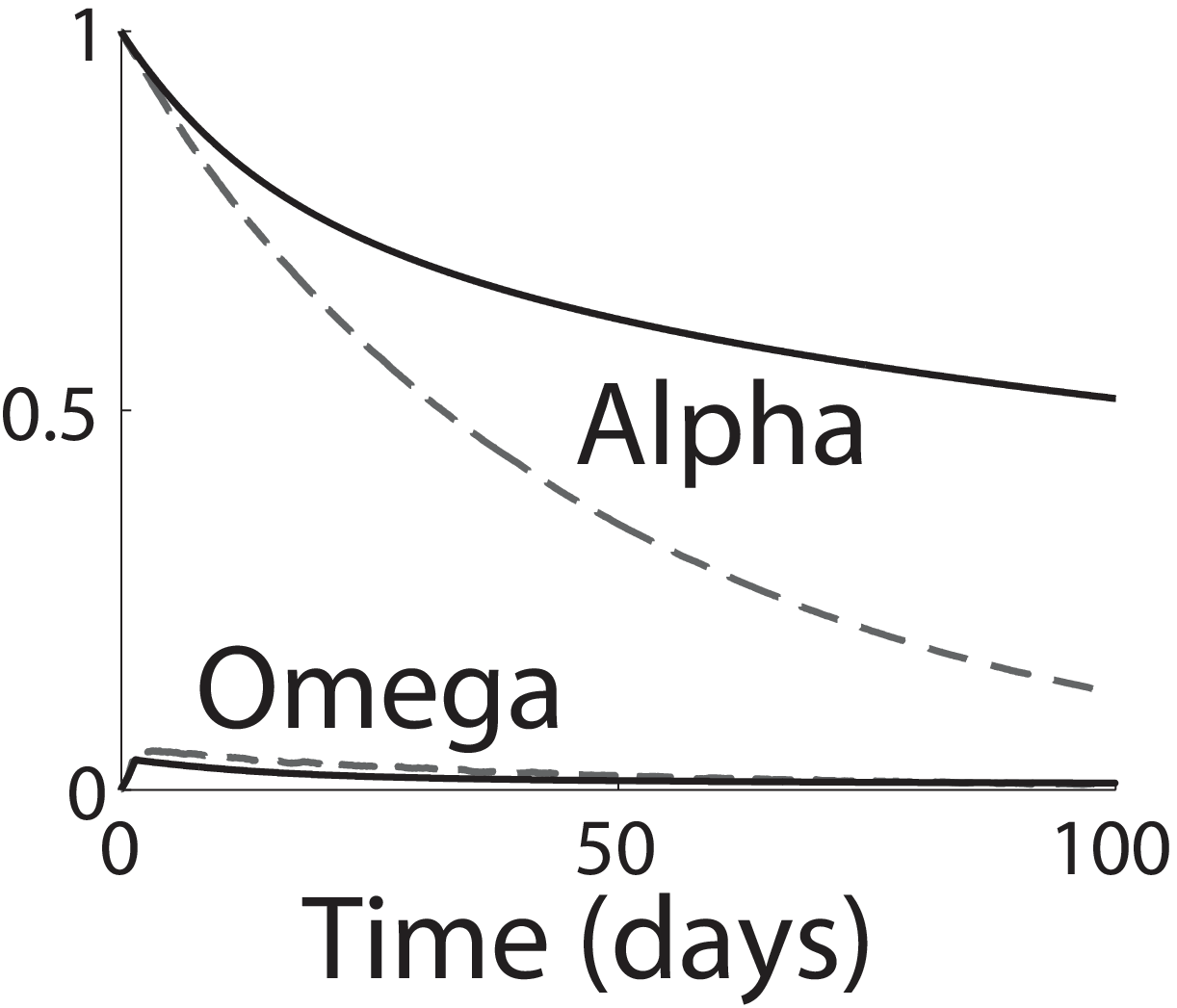} \\
{\bf (a)} & {\bf (b)} & {\bf (c)} \\
\end{tabular}
\end{center}
\caption{
Time plots of numerical solutions to Approximation 3 and simulations of the ABM.
Solutions to Approximation 3 are shown by solid black lines, and
simulations of the ABM are shown by dashed gray lines.
(a) Time plots for $d = 1.02$.
(b) Time plots for $d = 1.05$.
(c) Time plots for $d = 1.2$.
}
\label{figure:ComparisonApprox3toABM}
\end{figure}

Using the same initial conditions and parameters as before,
Figure~\ref{figure:ComparisonApprox3toABM} compares Approximation 3 and the ABM.
Although both Approximations 1 and 2 followed the behavior of the ABM very well, in this case, we see that combining the two approximations leads to quite different behavior.
%%%%-----------------------------------------------------------------------------------

\begin{figure}[htbp]
\begin{center}
\begin{tabular}{ccc}
\includegraphics[scale = 0.39]{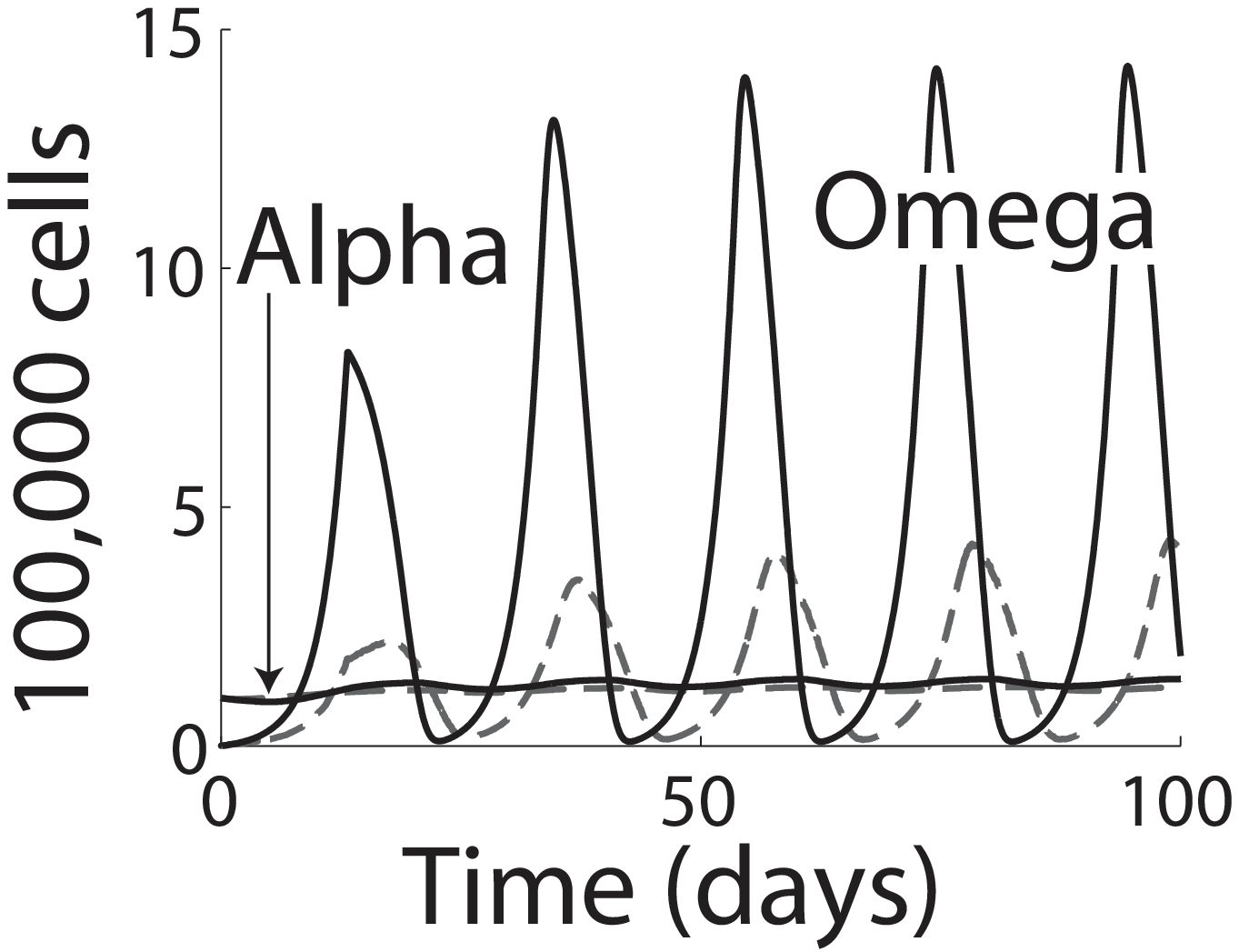} &
    \includegraphics[scale = 0.39]{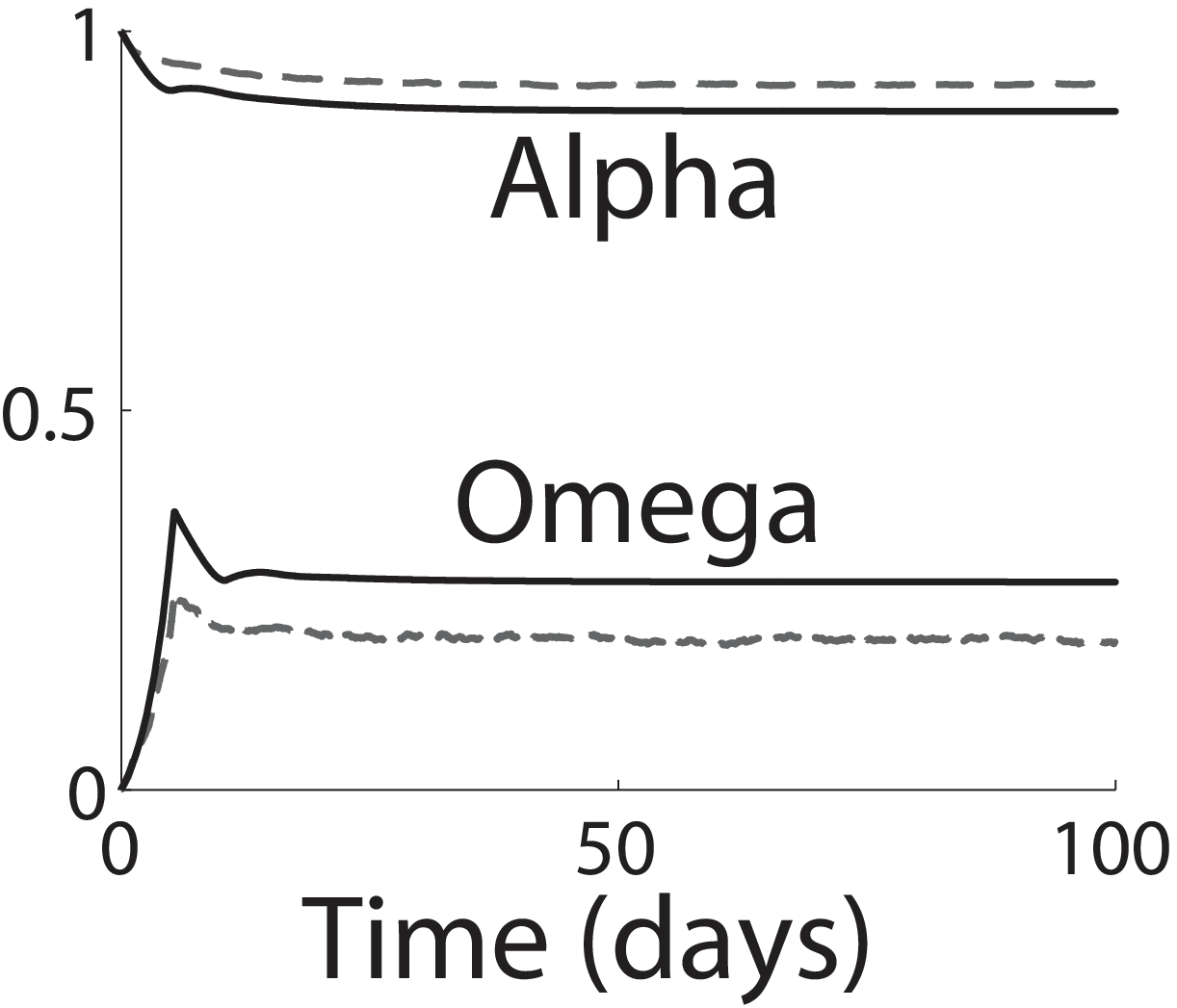} &
    \includegraphics[scale = 0.39]{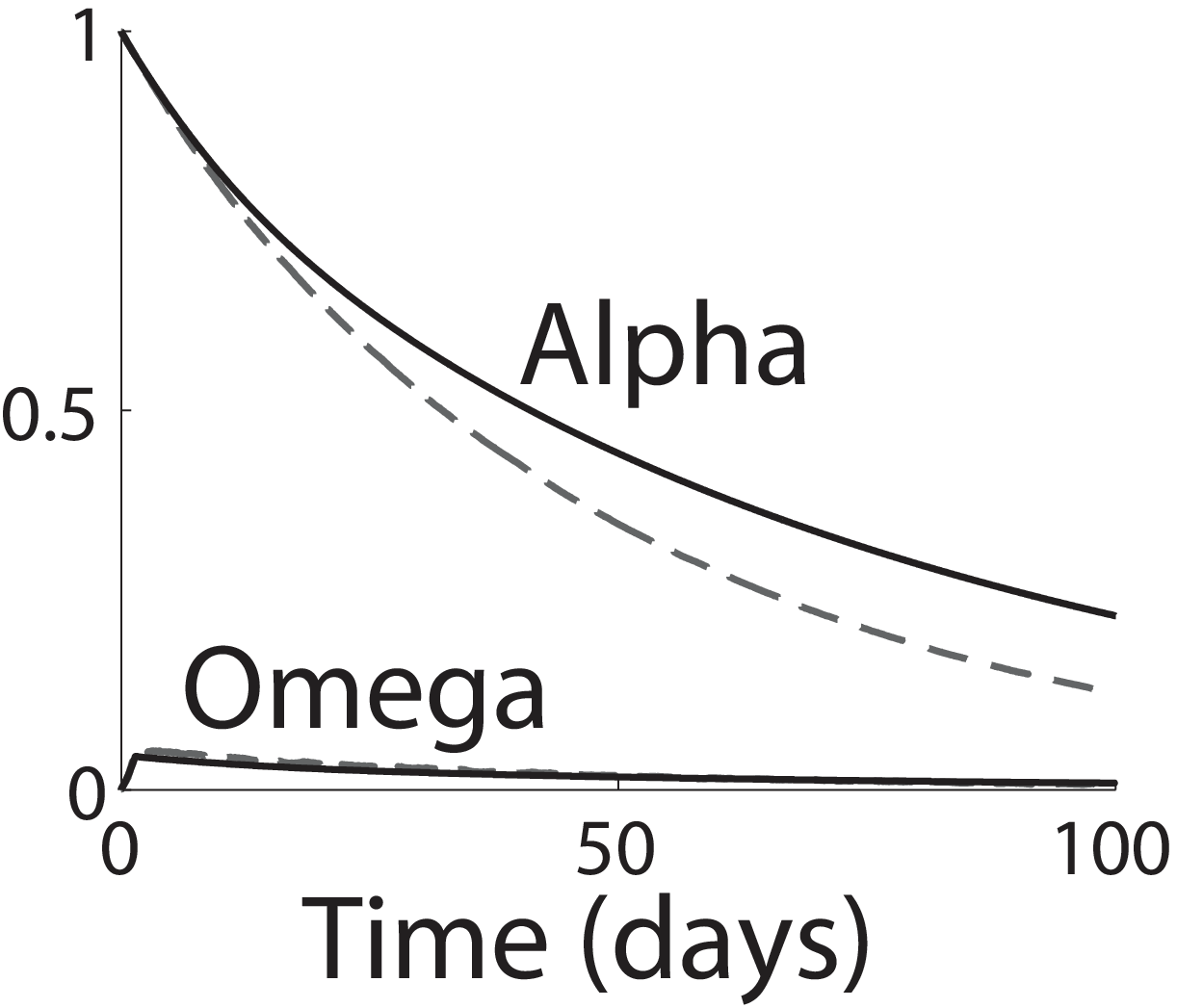} \\
{\bf (a)} & {\bf (b)} & {\bf (c)} \\
\end{tabular}
\end{center}
\caption{
Time plots of numerical solutions to Approximation 4 and simulations of the ABM.
Solutions to Approximation 4 are shown by solid black lines, and
simulations of the ABM are shown by dashed gray lines.
(a) Time plots for $d = 1.02$.
(b) Time plots for $d = 1.05$.
(c) Time plots for $d = 1.2$.
}
\label{figure:ComparisonApprox4toABM}
\end{figure}

Using the same initial conditions and parameters as before,
Figure~\ref{figure:ComparisonApprox4toABM} compares Approximation 4 and the ABM.
Approximation 4 behaves a lot like Approximation 3, so the added assumption that the transfer function $\alpha$ does not depend on $x$ does not affect the dynamics of the model very much.
%%%%-----------------------------------------------------------------------------------
%%%%-------------------------------------------------------------------------------------

\section{Analysis of Approximations 3 and 4}
\label{section:Analysis}
\setcounter{equation}{0}
\setcounter{figure}{0}
\setcounter{table}{0}
%%%%-------------------------------------------------------------------------------------
%%%%-----------------------------------------------------------------------------------
In this section, we study both the system \eqref{equation:Approx3alphastar}--\eqref{equation:Approx3omegastarBC} corresponding to Approximation 3 and the simpler system \eqref{equation:Approx4alphastar}--\eqref{equation:Approx4omegastarBC} corresponding to Approximation 4.
To conduct the analysis we use duality arguments related to the ''General Relative Entropy`` method introduced in \cite{BP,MMP} and used widely (see \cite{Calvez,Calvez2,CBBP1,CBBP2,MD} for other examples of applications). 
This method requires us to handle the eigenvalue problem and its adjoint, which we do first, and then use them to build entropy functionals.
Another method is to reduce the system to a delay differential equation (DDE), and we also comment on how this is possible.
%%%%-----------------------------------------------------------------------------------
%%%%-----------------------------------------------------------------------------------
\subsection{Link with a Delay Differential Equation}
\label{section:LinkWithDDE}
%%%%-----------------------------------------------------------------------------------
%%%%-----------------------------------------------------------------------------------
In the system (\ref{equation:Approx4alphastar})--(\ref{equation:Approx4omegastarBC}), we can interpret $\Omega^*$ as the maturing and proliferating stem cell compartment, whereas $A^*$ represents a reservoir of quiescent, completely immature stem cells. There are exchanges between these compartments, with the same rules as before, except that once a cell stops the maturation process and enters the quiescent compartment, it immediately becomes fully immature, but if it reenters the maturing and proliferating compartment, it must go through the entire maturation process again.
The main idea remains, however, unchanged, i.e., maturation is considered an invertible process as long as full maturity (at $x=1$) is not reached.

Although this assumption is original, and up to our knowledge has first been proposed in the ABM model of \cite{roeder2006}, the simplified system (\ref{equation:Approx4alphastar})--(\ref{equation:Approx4omegastarBC}), has some resemblance with the CML models based on delay differential equations (DDEs), proposed by Mackey \emph{et al.} \cite{adimy2005, bernard2003, bernard2004, colijn2005A, colijn2005B, pujo2004}.
Indeed, we can convert the PDE system (\ref{equation:Approx4alphastar})--(\ref{equation:Approx4omegastarBC}) into an equivalent DDE system by
applying the method of characteristics to (\ref{equation:Approx4omegastar}) to obtain
\begin{align*}
\Omega^*(x,t) & = \Omega^*\left( 0, t-\frac{x}{\rho_d} \right) \exp \left\{ \frac{b}{\rho_d}x
-\int_{t-x/\rho_d}^t \alpha(A^*(u)) du \right\} \\
& =\frac{\omega(\Omegatotal)}{\rho_d} A^* \left( t-\frac{x}{\rho_d} \right) \exp \left\{ \frac{b}{\rho_d}x -\int_{t-x/\rho_d}^t \alpha(A^*(u)) du \right\} ,
\end{align*}
for $1<\rho_d t$,
 a restriction that discards the initial data and simplifies the setting.
Integrating (\ref{equation:Approx4omegastar}) with respect to $x$ and using the equation above to express $\Omega(t,1)$ and $\Omega(t,0)$ in terms of $A^*,$ we obtain
\begin{align*}
\frac{d\Omegatotal}{d t} + \omega(\Omegatotal) e^{\f{b}{\rho_d} - C(t)} A^*(t - 1/\rho_d) - \omega(\Omegatotal) A^* 
& = \left( b - \alpha(A^*)  \right) \Omegatotal,
\end{align*}
where
$$
C(t) = \int_{t-\frac{1}{\rho_d}}^t \alpha(A^*(u)) du.
$$
Thus, we obtain the following DDE system:
\begin{align}
\label{eq:delay:1}
\frac{d\Omegatotal}{d t} 
& = \omega(t) A^*(t) -\omega(t) e^{\f{b}{\rho_d} - C(t)} A^*(t - 1/\rho_d) + \left( b - \alpha(t) \right) \Omegatotal(t),
    \\
\label{eq:delay:2}
\frac{dA^*}{dt}
& = - \omega(t) A^*(t) + \alpha(t) \Omegatotal(t), \\
\label{eq:delay:3}
\frac{dC}{dt}
& = \alpha(t) - \alpha(t-1/\rho_d),
\end{align}
where we use $\alpha(t)$ and $\omega(t)$ to denote $\alpha(A^*(t))$ and $\omega(\Omegatotal(t))$.
The characteristic equation obtained from the DDE system \eqref{eq:delay:1}--\eqref{eq:delay:3} coincides with the eigenvalue equation \eqref{eq:lambda} given by our analysis of local stability around the nonzero steady state of the PDE system \eqref{equation:Approx4alphastar}--\eqref{equation:Approx4omegastarBC}.
(See Section~\ref{section:stabilityanalysisnontrivialstate}.)
%%%%-----------------------------------------------------------------------------------
%%%%-----------------------------------------------------------------------------------
%-------------------------------------------------------------------------------------
%\subsection{Asymptotic Analysis}
%\label{section:AsymptoticAnalysis}
%MD CHANGE 03/15: title ''of Approximation 4``+ sentence of introduction
%-------------------------------------------------------------------------------------
\subsection{Existence of Steady States for Approximations 3 and 4}
%%%%-----------------------------------------------------------------------------------
The model based on Approximation 4 is simpler and we analyze it first. The method is extended to Approximation 3 afterwards.
\emph{Steady states} refer to time-independent solutions $(A^*,\Omega^*)$ of system \eqref{equation:Approx4alphastar}--\eqref{equation:Approx4omegastarBC}. They can be found and classified according to the model parameters.
\begin{proposition}\label{prop:steady}
Let $\alpha(\cdot)$ and $\omega(\cdot)$ be continuous positive bounded functions decreasing to zero.
Then zero is always a steady state, and
%\\
%$\bullet$ For $b=\rho_d,$  there is a unique nonzero steady state \emph{iff} $\alpha(A^*=0)>b$.  If $\alpha(A^*=0) %\leq b$, then zero is the unique steady state.
%\\
%$\bullet$ For $b\neq \rho_d,$
 there is a nonzero steady state (which is unique)  \emph{iff} $\alpha(0)>b^*$, where $b^*$ is uniquely determined by $b=b^*$ for $b=\rho_d$ and by the relation
\be \label{eq:defb*}
be^{-\f{b}{\rho_d}}=b^* e^{-\f{b^*}{\rho_d}},\qquad b\neq b^*, \qquad \text{for }\; b \neq \rho_d.
\ee
If $\alpha(0) \leq b^*$, then zero is the unique steady state.
\end{proposition}
\begin{figure}[ht]
\begin{center}
%\includegraphics[width=10cm, height=5cm]{Fig_F}
%\vspace{-1cm}
\includegraphics[scale=0.3]{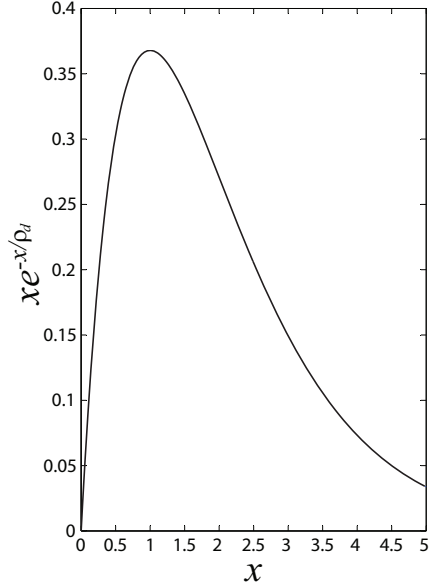}
\caption{\label{fig:F} The function $F(x)=xe^{-\f{x}{\rho_d}}$ used in (\ref{eq:defb*}). Its maximum is attained for $x=\rho_d$.}
\end{center}
\end{figure} 
\proof Obviously $ (\Omega=0, A^*=0)$ is a steady state.  We now consider possible nonzero steady states. Solving (\ref{equation:Approx4omegastar}), we find that nonzero steady state solutions $(\tilde\Omega(x),\tilde A)$ satisfy
\be
\label{eq:Omegatilde}
\tilde \Omega (x) = \tilde \Omega (0) e^{\f{b-\alpha(\tilde A)}{\rho_d}x},
\ee 
where $\Omega (0) \neq 0$. Substituting (\ref{eq:Omegatilde}) into (\ref{equation:Approx4alphastar}) and (\ref{equation:Approx4omegastarBC}), we determine the solution by the values $(\tilde A,\Omegatotal )$ that satisfy
\be\label{eq:steadystate}
\tilde A= \frac{\alpha(\tilde A) \Omegatotal}{\omega(\Omegatotal)}, \qquad  \rho_d  =\alpha(\tilde A) \int_0^1 e^{\f{b-\alpha(\tilde A)}{\rho_d}x}dx, \qquad \Omegatotal = \tilde\Omega (0) \int_0^1 e^{\f{b-\alpha(\tilde A)}{\rho_d}x}dx.
\ee
Thus, we have  two cases:

\noindent\textbf{Case 1.} $b=\rho_d$. This is equivalent to  $b = \alpha(\tilde A)$ by the second equation of (\ref{eq:steadystate}). Therefore $b^*=b= \alpha(\tilde A)$ and we can choose an $A^* = \tilde{A}\neq 0$ under the condition $\alpha(0)>b^*$. Then the monotonicity condition on $\omega(\cdot)$ allows us to find a unique $\Omegatotal$ by the first relation in (\ref{eq:steadystate}).

\noindent\textbf{Case 2.}  $b\neq\rho_d$ and thus  $b \neq \alpha(\tilde A)$. Since the second relation also implies 
$$
be^{-\f{b}{\rho_d}} = \alpha(\tilde A) e^{-\f{\alpha(\tilde A)}{\rho_d}}.
$$
we clearly have $b^*= \alpha(\tilde A)$  and the result follows as before. 
\qed
\begin{proposition}\label{prop:steady3}
Let $\alpha(x, \cdot)$ and $\omega(\cdot)$ be continuous positive bounded functions decreasing to zero (for all $x$ in case of $\alpha$).
Then zero is always a steady state, and there is a nonzero steady state (which is unique) \emph{iff} 
\be \label{eq:condstst3 }
\int_0^1 \alpha(x,0)e^{\int\limits_0^x \frac{b-\alpha(y,0)}{\rho_d}dy}dx> \rho_d \;  \Longleftrightarrow \; \int_0^1 e^{\int\limits_x^1 \frac{\alpha(y,0)-b}{\rho_d}dy} dx> \frac{\rho_d}{b}.  
\ee
\end{proposition}
We leave to the reader that this condition is equivalent to $\alpha(y,0) >b^*$ in the case $\alpha$ is independent of $x$ as stated in Proposition \ref{prop:steady}.

\proof The explicit solution is now given by 
$$\Omega^*(x)= \Omega^*(0) e^{\int\limits_0^x \frac{b-\alpha(y,A^*)}{\rho_d}dy}.$$
Therefore,  as in the proof of Proposition \ref{prop:steady}, one reduces the claim to first finding $A^*$ that solves
$$
1= \int_0^1 \frac{\alpha(x,A^*)}{\rho_d}e^{u(x,A^*)} dx,
$$
with $u(x,A^*)=\int_0^x \frac{b-\alpha(y,A^*)}{\rho_d}dy$.  Because $\frac{\alpha(x,A^*)}{\rho_d}=\frac{b}{\rho_d}- u_x'(x,A^*)$, integration by parts reduces our condition to 
$$
1= \frac{b}{\rho_d} \int_0^1e^{u(x,A^*)}dx+e^{u(0,A^*)} - e^{u(1,A^*)},
$$
or also
$$
 \frac{\rho_d} {b}= \int_0^1e^{u(x,A^*)-u(1,A^*)}dx= \int_0^1e^{\int\limits_x^1  \frac{\alpha(y,A^*)-b}{\rho_d}dy}dx=:H(A^*).
$$
This function $H(\cdot)$ is decreasing and satisfies
$$
H(\infty)=  \frac{\rho_d} {b} (1-e^{ -\rho_d/b})< \frac{\rho_d} {b}.
$$
Therefore our condition (\ref{eq:condstst3 }) on $H(0)$ is necessary and sufficient for the existence of a solution  $A^*>0$. The end of the proof is the same as before.
\qed
%%%%-------------------------------------------------------------------------------------

\subsection{An eigenvalue problem for Approximation 3}
\label{sec:eigen}
Our analysis of {\em a priori} bounds and stability uses the eigenelements of the linear problem related to various states. This section introduces the relevant material. 

Here we fix $A^*$ and $\Omega^*,$ we set $\alpha(x) = \alpha(x, A^*)$ and $\omega = \omega(\Omegatotal),$ and we consider the eigenvalue problem to find $(\lambda\in \R, \Omega(x), A)$ that satisfy 
\begin{align}
\label{eq:eigen:1}
& \rho_d \frac{\partial \Omega}{\partial x} = 
    \left( b -\lambda- {\alpha}  \right) \Omega, \qquad 0\leq x\leq 1,
    \\
\label{eq:eigen:2}
& 0 = \int_0^1 {\alpha} \Omega (x) dx -(\lambda+ {\omega} ) A, 
\\ 
\label{eq:eigen:3}
& \Omega(0) = \frac{{\omega} } {\rho_d} A,\qquad \int_0^1 \Omega(x)dx=1.
\end{align}
Notice that we can simplify the problem by combining the last two relations to 
\be\label{eq:eigen:S:2}
(\lb +\omega)\Omega(0)
 = \f{\omega}{\rho_d} \int_0^1 {\alpha} \Omega (x) dx,\qquad \int_0^1 \Omega(x)dx=1.
\ee

\begin{proposition}\label{prop:eigen}
Let $\alpha(x)>0$ and $\omega>0$. There exists a unique solution $\lambda$, $\Omega(x)>0$, $A>0$ to (\ref{eq:eigen:1})--(\ref{eq:eigen:3}).  This $\lambda$ is the unique real eigenvalue, and it is positive \emph{iff} 
\begin{equation}
 \label{eq:lbpos}
\int_0^1 \alpha(x) e^{\f{1}{\rho_d}\int\limits_0^x (b-\alpha(s))ds}dx >\rho_d,
\end{equation}
and it is negative \emph{iff}
\begin{equation}
 \label{eq:lbneg}
\int_0^1 \alpha(x) e^{\f{1}{\rho_d}\int\limits_0^x (b-\alpha(s))ds}dx <\rho_d.
\end{equation}
If $\alpha(x)=\alpha$ is a constant, then (\ref{eq:lbpos}) reduces to $\alpha>b^*$ 
where $b^*$ is defined by (\ref{eq:defb*}). The eigenvalue equals 0 \emph{iff} $\alpha=b^*.$
\end{proposition}
\proof 
One calculates explicitly from (\ref{eq:eigen:1}) that
$$
\Omega(x)=\Omega(0) e^{\f{1}{\rho_d} \int\limits_0^x (b-\alpha(s)-\lb)ds},
$$
and together with (\ref{eq:eigen:S:2}), that the eigenvalue $\lb$ must satisfy
$$\rho_d \left( \f{\lb}{\omega}+1 \right) = \int_0^1 \alpha(x) e^{\f{1}{\rho_d} \int\limits_0^x (b-\alpha(s)-\lb) ds} dx.$$
We can write this relation in the form $G(\lb)=L(\lb),$ with 
$$
G(\lb)=\rho_d \left( \f{\lb}{\omega}+1 \right) ,\qquad L(\lb)=\int_0^1 \alpha(x) e^{\f{1}{\rho_d} \int\limits_0^x (b-\alpha(s)-\lb) ds} dx.
$$
The function $G$ increases linearly from $-\infty$ to $+\infty,$ whereas $L(\lb)$ decreases continuously from $+\infty$ to $0$.  Hence, the two curves intersect at a unique value $\lb\in \R.$ Moreover, $\lb>0$ \emph{iff} $G(0)<L(0),$ which gives Equation (\ref{eq:lbpos}).
If $\alpha(x)$ is constant, we obtain
$$
1 < \f{\alpha}{b-\alpha}\left( e^{\f{b-\alpha}{\rho_d}}-1 \right),
$$
which, with the notation $F(z)=ze^{-z/\rho_d}$,  simplifies to
$$
\frac{1}{b-\alpha} F(b) < \frac{1}{b-\alpha} F(\alpha)
$$
If $b>\alpha$, this inequality is equivalent to $F(\alpha) > F(b)$, so $\alpha > \min(b,b^*) = b^*.$ If $b<\alpha,$ the inequality is equivalent to $F(\alpha)<F(b)$, so $\alpha > b^*$.  Thus, in both cases, $\alpha > b^*$.
\qed 

\noindent {\bf Remark.} The adjoint problem of (\ref{eq:eigen:1})--(\ref{eq:eigen:3}) is
\begin{align}\label{eq:eigen:ad:1}
- \rho_d \frac{\partial \phi}{\partial x} & = 
    \left( b -\lambda- {\alpha} (x) \right) \phi + \alpha(x)\Psi,
    \\
\label{eq:eigen:ad:2}
(\lambda + \omega)\Psi&=\omega \phi(0),
\\  \label{eq:eigen:ad:3}
\phi(1) & = 0,
\end{align}
and its solution $\phi$ can be explicitly calculated as follows:
\be \label{eq:phi} \phi(x)=\phi(0)e^{-\f{1}{\rho_d}\int\limits_0^x(b-\alpha(s)-\lb)ds} \left(
1- \f{\int_0^x \alpha(s) e^{\f{1}{\rho_d}\int\limits_0^s (b-\alpha(\sigma)-\lb)d\sigma}ds}{\int_0^1 \alpha(s) e^{\f{1}{\rho_d}\int\limits_0^s (b-\alpha(\sigma)-\lb)d\sigma}ds}\right).
\ee 
%%%%-----------------------------------------------------------------------------------
%%%%-----------------------------------------------------------------------------------
%%%%-----------------------------------------------------------------------------------
%-------------------------------------------------------------------------------------
\subsection{Uniform Bounds}
%%%%-----------------------------------------------------------------------------------
%%%%-----------------------------------------------------------------------------------
The aim of this section is to show that the population remains bounded under the previous assumption that $\alpha$ vanishes at infinity (this can be relaxed as shown below). This relies  on the eigenvalue problem and the property that the eigenvalue is negative for large populations.
%%%%-----------------------------------------------------------------------------------
%%%%-----------------------------------------------------------------------------------
\begin{proposition}\label{prop:bound}
If $\omega(\cdot)$ is bounded, $\alpha(A,x)$ is a nonincreasing function of $A$ and for $A$ large enough
\be\label{as:infty}
 \int_0^1 \alpha(x,A)e^{\int\limits_0^x \f{b-\alpha(y,A)}{\rho_d}dy} dx <\rho_d, 
\ee
then any solution $\bigl(\Omega^*(t,x),A^*(t)\bigr)$ of (\ref{equation:Approx3alphastar})--(\ref{equation:Approx3omegastarBC}) remains bounded for all $t\geq 0$, i.e., $A^*\in  L^\infty(0,\infty)$ and  $\Omega^*\in L^\infty((0,\infty)\times (0,1))$.
\end{proposition}
%%%%-----------------------------------------------------------------------------------
\proof
Since (\ref{as:infty}) holds, we can choose $\tilde A$ large enough so that the condition (\ref{eq:lbneg}) is fulfilled for $\alpha(x)=\alpha(x,\tilde A)$.
According to Proposition \ref{prop:eigen}, there exists a unique eigenvalue $ \lb <0$ and a unique positive couple $(A, \Omega)$ that solves (\ref{eq:eigen:1})--(\ref{eq:eigen:3}) with $\alpha(x)=\alpha(x,\tilde A)$.
We also denote by $\phi$ the solution of the adjoint eigenproblem (\ref{eq:eigen:ad:1})--(\ref{eq:eigen:ad:3}). 
Testing the solutions $\Omega^*(t,x)$ of Problem (\ref{equation:Approx3alphastar})--(\ref{equation:Approx3omegastarBC}) against $\phi$ gives
\begin{equation*}\begin{array}{ll}
\f{d}{dt}\int  \phi \Omega^* dx &= \int \phi \bigl(b-\alpha(x,A^*)\bigr)\Omega^* dx+ \int \rho_d \Omega^* \f{\p}{\p x}  \phi dx + \rho_d \Omega^*(t,0) \phi (0) \\
& =\int \Omega^*\bigl(\lb\phi- \alpha (x,\tilde A)\Psi \bigr)dx + \int \bigl(\alpha(x,\tilde A)-\alpha(x,A^*(t))\bigr)\Omega^*\phi dx \\
& \quad + \rho_d \Omega^*(t,0) \phi(0)
\\
&= \int\biggl(\Omega^*\bigl(\lb \phi-\f{\alpha(x,\tilde A)\omega(\Omegatotal)}{\lb+\omega(\Omegatotal)}\phi(0)  \bigr) +  \bigl(\alpha(x,\tilde A)-\alpha(x,A^*(t))\bigr) \Omega^* \phi\biggr) dx \\
& \quad + \phi(0)\bigl(-\f{dA^*}{dt} + \int_0^1 \alpha(x,A^*(t)) \Omega^* dx \bigr)
\\
& = \int\Omega^*\bigl(\lb \phi+\lb \f{\alpha(x,\tilde A) \Psi}{\omega(\Omegatotal)}  \bigr) dx+ \int \bigl(\alpha(x,\tilde A)-\alpha(x,A^*(t))\bigr) \Omega^* \phi dx 
\\
& \quad + \phi(0) \bigl(-\f{dA^*}{dt} 
+\int \bigl(\alpha (x,A^*(t))-\alpha (x,\tilde A))\Omega^* dx)
\\
& =\int\Omega^*\bigl(\lb \phi +\lb \f{\alpha(x,\tilde A)\Psi}{\omega(\Omegatotal)} \bigr) dx+ \int  \bigl(\alpha(x,\tilde A)-\alpha(x,A^*(t))\bigr)\Omega^*  \bigl(\phi(x)-\phi(0)\bigr) dx \\
& \quad -\phi(0)\f{d}{dt} A^* 
\end{array}
\end{equation*}
Hence, we have obtained the following equality:
\begin{align*}
\begin{split}
 \f{d}{dt}\left(\int \phi \Omega^* dx + \phi(0) A^*\right) & = \lb \int\Omega^*\left( \phi + \f{\alpha(x,\tilde A)\Psi}{\omega(\Omegatotal)} \right) dx \\
& \quad + \int \Omega^*  \left(\phi(x)-\phi(0)\right) \left(\alpha(x,\tilde A)-\alpha(x,A^*(t))\right) dx.
\end{split}
\end{align*}
%MD CHANGE 03/15: PROOF BY Benoit PERTHAME!
We consider the quantity $S(t)=\int \phi (x) \Omega^*(x,t) dx + \phi(0) A^*(t) - \mu \min \left(\tilde A, A^*(t)\right),$ where $\mu$ is a constant such that 
\begin{equation}\label{eq:mu}\forall x\in [0,1],\qquad \mu \alpha(x,\tilde A) \geq (-\lb+1) \phi + ||\alpha(.,.)||_\infty ||\phi||_\infty.\end{equation}
Recalling that $\lb <0$, we obtain
\begin{align}
\begin{split}
\frac{dS}{dt} & \leq \lb \int \Omega^* \phi  
+
 \int \Omega^*  \left(\phi(x)-\phi(0)\right) \left(\alpha(x,\tilde A)-\alpha(x,A^*(t))\right) dx
 \\
& \quad + \mu {\mathbf{1}}_{A^*(t)\leq \tilde A} \left(\omega A^*(t) - \int_0^1 \alpha(x,A^*(t))\Omega^*(t,x) dx \right).
\label{eq:GRE2}
\end{split}
\end{align}
We set $S_0 = \phi(0) \tilde A + \mu ||\omega||_\infty \tilde A$ and we prove that $S(t)\leq \max(S(0),S_0).$  Indeed, If $S(t)\geq S_0,$ we have two cases:

\noindent $\bullet${\bf First case: $A^*(t) \geq \tilde A.$} The last term in the RHS of \eqref{eq:GRE2} vanishes, the first two terms are less than or equal to zero, so $\frac{dS}{dt}\leq 0$ and $S(t)$ decreases.

\noindent $\bullet$ {\bf Second case: $A^*(t) \leq \tilde A.$} Equation \eqref{eq:GRE2} implies (recall that $\alpha(x,A^*) \geq \alpha(x,\tilde A)$ )
$$ \frac{dS}{dt} \leq \lb \int \Omega^* \phi  
+
 \int \Omega^* \phi(0)\alpha(x,A^*(t))
+\mu\omega A^*(t)
-\mu\int_0^1 \alpha(x,A^*(t))\Omega^*(t,x) dx .
$$
We deduce, using \eqref{eq:mu}, that
$$ \frac{dS}{dt} \leq \mu ||\omega||_\infty \tilde A - \int_0^1 \Omega^* \phi dx\leq S_0 - S(t) \leq 0.$$

This shows that $S(t)$ decreases in both cases, which proves that $S(t) \leq \max(S(0),S_0)$ remains constantly bounded. As a consequence $\int \phi (x) \Omega^*(x,t) dx + \phi(0) A^*(t)$ is uniformly bounded. Then, since $A^*(t)$ is bounded, in  (\ref{equation:Approx3omegastarBC}) the boundary term $\Omega^*(0,t)$ is also uniformly bounded and thus $\Omega^*\in L^\infty((0,\infty)\times (0,1))$. 
\qed

%--------------------------------------------------------------------------------------
\subsection{Stability Analysis near Zero}

\begin{proposition}\label{prop:Zero}
If $\alpha(x,A)$ is a positive decreasing function of $A$, and if 
\begin{equation}\label{ineq:stab}
\int_0^1 \alpha(x,0) e^{\f{1}{\rho_d}\int\limits_0^x (b-\alpha(s,0))ds}dx \leq \rho_d,
 \end{equation}
then the zero steady state is globally attractive. 
On the contrary, if 
\begin{equation}\label{ineq:unstab}\int_0^1 \alpha(x,0) e^{\f{1}{\rho_d}\int_0^x (b-\alpha(s,0))ds}dx >\rho_d,
 \end{equation}
then the zero steady state is unstable.
\end{proposition}
\noindent {\bf Remark.} According to the proof of Proposition \ref{prop:steady3}, the stability condition (\ref{ineq:unstab}) is equivalent to condition (\ref{eq:condstst3 }) for the existence of a nonzero steady state. Therefore the zero steady state is stable exactly when it is the unique steady state. Also, the properties of $\omega(0,\Omegatotal)$ do not influence the stability of the zero steady state.

\proof  We consider the linear eigenvalue problem (\ref{eq:eigen:1})--(\ref{eq:eigen:3}) at zero and let $\phi$, given by (\ref{eq:phi}), be the solution to the adjoint eigenvalue problem.
Let $\left(\Omega(t,x), A(t)\right)$ denote a solution to the time-dependent problem (\ref{equation:Approx4alphastar})--(\ref{equation:Approx4omegastarBC}). As in the proof of Proposition \ref{prop:bound}, we test $\Omega$ against $\phi$  and define 
$v(t)=\int \phi(x) \Omega(x,t)  dx + \phi(0) A(t)$. We have 
%\begin{align*}
%\f{d}{dt}\int \phi \Omega dx &= \int \phi \bigl(b-\alpha(x,B)\bigr)\Omega dx+ \int \rho_d \Omega \f{\p}{\p x} \phi dx + \rho_d \Omega(0) \phi (0)
%\\
%&=\int \Omega\bigl(\lb\phi- \alpha (x,0)\Psi \bigr)dx + \int \bigl(\alpha(x,0)-\alpha(x,B)\bigr)\Omega \phi dx + \rho_d \Omega(0)\phi(0)
%\\
%&= \int\biggl(\Omega\bigl(\lb \phi-\f{\alpha(x,0)\omega(0)}{\lb+\omega(0)}\phi(0)  \bigr) +  \bigl(\alpha(x,0)-\alpha(x,B)\bigr) \Omega \phi\biggr) dx + \phi(0)\bigl(-\f{dB}{dt} + \int_0^1 \alpha(B,x) \Omega dx \bigr)
%\\
%&= \int\biggl(\Omega\bigl(\lb \phi-\f{\alpha(x,0)\omega(0)}{\lb+\omega(0)}\phi(0)  \bigr) +  \bigl(\alpha(x,0)-\alpha(x,B)\bigr) \Omega \phi\biggr) dx + \phi(0)\bigl(-\f{dB}{dt} + \int_0^1 \alpha(B,x) \Omega dx \bigr)
%\\
%&=  \int\Omega\bigl(\lb \phi+\lb \f{\Psi}{\omega(0)}  \bigr) dx+ \int \bigl(\alpha(x,0)-\alpha(x,B)\bigr) \Omega \phi dx 
%\\
%&+\phi(0) \bigl(-\f{dB}{dt} 
%+\int \bigl(\alpha (x,B)-\alpha (x,0))\Omega dx)
%\\
%&=\int\Omega\bigl(\lb \phi +\lb \f{\Psi}{\omega(0)} \bigr) dx+ \int  \bigl(\alpha(x,0)-\alpha(x,B)\bigr)\Omega  \bigl(\phi(x)-\phi(0)\bigr) dx -\phi(0)\f{d}{dt} B 
%\end{align*}
\begin{align}
\begin{split}\label{eq:GRE:0}
 \f{d}{dt}v(t)& = \lb \int\Omega\left( \phi + \f{\alpha(x,0)\Psi}{\omega(0)} \right) dx \\
& \quad 
+ \int \Omega \left(\phi(x)-\phi(0)\right) \left(\alpha(x,0)-\alpha(x,A)\right) dx.
\end{split}
\end{align}
If $\lb<0$, which is equivalent to (\ref{ineq:stab}) according to Proposition \ref{prop:eigen},  \eqref{eq:phi} shows that  $\phi$ is decreasing, and since we have supposed that $\alpha$ decreases with respect to $A$, we can conclude that
$$
\f{d}{dt}\int \phi \Omega dx + \phi(0) A \leq \lb \int \phi \Omega dx.
$$
 It follows that $\int \phi \Omega dx + \phi(0) A$ is decreasing.  Since this integral is nonnegative, it tends to a limit, and at infinity, $ \int \phi \Omega dx$ tends to $0$. Therefore, from the uniform bounds,  $\Omegatotal(t) \underset{t \to \infty}{\longrightarrow}  0$, and  thus $A(t) \underset{t \to \infty}{\longrightarrow}   0$.
Thus, the zero steady state is globally attractive.

If $\lb=0,$ we 
still have that $v(t)$ is nonincreasing, and one has
$$
\f{d}{dt}v(t)  = \int \Omega(x,t) \left( \phi(x)-\phi(0) \right)(\alpha(x,0)-\alpha(x,A))dx.
$$
Since $\phi(x)-\phi(0) <0$ and $\alpha(x,0)-\alpha(x,A) >0$ for $x<1$ and $A>0,$ it implies that either
$\Omega (x<1,+\infty) = 0$ or $A(+\infty) = 0$.  We conclude using 
(\ref{equation:Approx4alphastar})--(\ref{equation:Approx4omegastarBC}) that $\Omega = A = 0$. 

If $\lb>0,$ 
%let us suppose that there exists a sequence $t_n$ such that 
%$S(t_n)=\int \phi\Omega dx + \phi(0) A (t_n) \leq \f{1}{n}.$ Equation (\ref{eq:GRE:0}) and Assumption (\ref{as:alpha}) give, as %soon as $t_n$ is large enough so that $A(t_n)\leq A_0,$ and denoting $\alpha_m = \min_{x\in[0,1]} \alpha(x,0) >0$
%$$\f{d}{dt} S_n \geq  \int \Omega(t_n,x) dx \biggl(\lb \phi + \lb \f{\Psi\alpha_m}{\omega(0)} -f(A(t_n)) \phi(0) \biggr).$$
%For $A(t_n)$ small enough, since $f$ vanishes, it implies either $\int \Omega(t_n,x) dx =0,$ what is absurd (cf. the explicit %formula given by the method of characteristics) as soon as $\Omega^{in} \neq 0,$ or 
%$$\f{d}{dt} S_n >0.$$
%Hence, it is possible to find another sequence $t_n'$ with $t_n'<t_n,$ $S(t_n')<S(t_n)\leq \f{1}{n},$ and $\f{d}{dt} S(t'_n)=0.$ The previous calculation shows that such a sequence cannot exist; hence $S(t)$ is bounded from below, 
we define the constant $C=\|\phi'\|_\infty \sup |\alpha'_A(x,A)|/\phi(0)$ where the $sup$ is taken on $0\leq x \leq 1$ and $0\leq  \phi(0) A \leq \infty$. Then, we write 
\begin{align*}
\begin{split}
 \f{d}{dt}v(t) &\geq \int_0^1 \phi \Omega  [ \lb -C \phi(0)A(t)]dx
 \\
 & \geq \int_0^1 \phi \Omega  [ \lb -Cv(t)]dx.
\end{split}
\end{align*}
This proves that $v(t)$ is increasing whenever $v(t) \leq \lb/C$
and the zero steady state is thus unstable.
\qed

{\bf Remark.} We have not made any assumption on $\omega(\Omegatotal)$.  Indeed, the population in the $A^*$ compartment neither divides nor dies. Hence, whether cells in $A^*$ transfer quickly or slowly into the $\Omega$ compartment does not change the behavior at infinity or at zero.
On the contrary, we can expect it influences the stability of the nonzero steady state, whenever it  exists.

%%%%--------------------------------------------------------------------------------------

%MD CHANGE 15/03: title ''for Approximation 4``
\subsection{Stability analysis near the nonzero steady state}
\label{section:stabilityanalysisnontrivialstate}

To study the local stability of the operator around the nontrivial steady state, denoted $\tilde{\Omega},\tilde A,$ we consider a small perturbation $\delta\Omega=\Omega - \tilde\Omega,$ $\delta A = A - \tilde A$, where $A(t,x)$ and $\Omega(t,x)$ are solutions of System \eqref{equation:Approx4alphastar}--\eqref{equation:Approx4omegastarBC}.

For the sake of simplicity, we denote $\overline{\tilde \Omega}=\int \tilde\Omega dx,$ $\alpha=\alpha \big(\tilde A \big),$ $ \alpha'=\frac{d\alpha}{dA} \big(\tilde A \big),$ $\omega = \omega\big(\overline{\tilde{\Omega}} \big)$ and $\omega'=\frac{d\omega}{d\Omegatotal} \big(\overline{\tilde\Omega} \big)$.

At first order in $\delta A$ and $\delta \Omega$, the perturbation $\delta \Omega, \;\delta A$ satisfies the following equation:
\begin{align}\label{eq:delta}
\frac{\partial \delta\Omega}{\partial t} + \rho_d \frac{\partial \delta\Omega}{\partial x} & = 
    \left( b - {\alpha}\right) \delta\Omega -{\alpha}' \tilde \Omega \delta A ,
    \\
\begin{split}
\frac{d\delta A}{dt}
& = ({\alpha}-\omega' \tilde A)\int_0^1  \delta \Omega (t,x) dx
+ (\int_0^1 \alpha'   \tilde \Omega (x) dx -\omega)\delta A,
\end{split} \\ \label{eq:delta3}
\delta\Omega(0,t) & = \frac{{\omega}}{\rho_d} \delta A(t)
+ \frac{{\omega}'\tilde A} {\rho_d} \int_0^1 \delta \Omega(t,x) dx .
\end{align}
We now consider the eigenvalue problem related to the linear system \eqref{eq:delta}--\eqref{eq:delta3}.
Using $\alpha(\tilde A)=b^*={b}e^{-\f{b-b^*}{\rho_d}}$ (cf. Equation \eqref{eq:defb*}) and the definition of $\tilde \Omega$ in \eqref{eq:Omegatilde}, we calculate
$$
\int_0^1 \tilde \Omega (x) dx= \tilde \Omega (0) \rho_d \f{e^{\f{b-b^*}{\rho_d}}-1}{b-b^*}
=\f{\tilde\Omega (0) \rho_d}{b^*}.
$$
We also recall that $\tilde \Omega (0) \rho_d = \tilde A \omega$ by \eqref{eq:steadystate}.  Hence, the eigenvalue problem writes 
\begin{align}\label{eq:delta:eigen}
\lb \delta\Omega + \rho_d \frac{\partial \delta\Omega}{\partial x} & = 
    \left( b - b^* \right) \delta\Omega -{\alpha}'\delta A \tilde \Omega (0)e^{\f{b-b^*}{\rho_d}x},
    \\ \label{eq:delta:eigen:2}
\lb\delta A
& = b^* \int_0^1  \delta \Omega (x) dx 
+   \alpha'\delta A\f{\tilde\Omega (0) \rho_d}{b^*} 
- \f{\tilde \Omega(0) \rho_d}{\tilde A}  \delta A 
- \omega' \tilde A\int_0^1 \delta \Omega(x) dx , 
\\ \label{eq:delta:eigen:3}
\delta\Omega(0) & = \frac{\tilde \Omega (0)} {\tilde A} \delta A
+ \frac{{\omega}'\tilde A} {\rho_d} \int_0^1 \delta \Omega(x) dx.
\end{align}
(We still denote the eigenvector by $(\delta A,\;\delta \Omega)$.)
By explicitly solving \eqref{eq:delta:eigen}, we obtain
$$\delta \Omega (x) =e^{\f{b-b^*-\lb}{\rho_d}x} \left(\delta \Omega(0)  - \alpha' \tilde\Omega (0)\delta A \f{e^{\f{\lb}{\rho_d}x} -1}{\lb}\right).$$
We impose the restriction
\be\label{eq:int1}
\int_0^1 \delta\Omega (x) dx=1= 
\f{\delta\Omega (0) \rho_d }{b-b^*-\lb}(\f{b}{b^*}e^{-\f{\lb}{\rho_d}}-1)
 - \f{\alpha' \tilde\Omega (0)\delta A \rho_d}{\lb}\f{1}{b^*}\f{b(1-e^{-\f{\lb}{\rho_d}})-\lb}{b-b^*-\lb},
\ee 
and obtain from \eqref{eq:delta:eigen:2} that
$$\lb \delta A = b^* + \alpha' \delta A  \f{\tilde \Omega (0)}{b^*} \rho_d 
- \f{\tilde\Omega(0)}{\tilde A} \rho_d \delta A - \omega' \tilde A.$$
%MD CHANGE 23/04: erase the following detail:
%Hence,
%$$\delta A= \f{b^* -\omega' \tilde A}{\lb - \alpha' \f{\tilde\Omega(0)}{b^*} \rho_d + \f{\tilde \Omega (0)} {\tilde A} \rho_d }.$$
In addition, \eqref{eq:delta:eigen:3} yields
$$ \delta \Omega(0)= \f{\omega' \tilde A}{\rho_d} + \f{\tilde \Omega (0)}{\tilde A} \f{b^* -\omega' \tilde A} {\lb - \alpha' \f{\tilde \Omega (0)}{b^*} \rho_d + \f{\tilde\Omega (0)}{\tilde A} \rho_d},
$$
which we substitute into \eqref{eq:int1} to obtain
\begin{align}
\begin{split}
\label{eq:lambda}
1=f(\lb) & := \f{\rho_d}{b-b^* -\lb} \left(\f{b}{b^*}e^{-\f{\lb}{\rho_d}}-1\right)
\left(\f{\omega' \tilde A}{\rho_d} + \f{\tilde \Omega (0)}{\tilde A} \f{b^* - \omega'\tilde A}{\lb-\f{\alpha' \tilde \Omega (0)}{b^*} \rho_d + \f{\tilde \Omega (0)}{\tilde A}\rho_d}\right)
\\
& \quad - \f{\alpha' \tilde\Omega(0) \rho_d}{b^* \lb} \f{b \left( 1-e^{-\f{\lb}{\rho_d}} \right)-\lb}{b-b^*-\lb}
\f{b^* -\omega'\tilde A}{\lb - \alpha' \f{\tilde \Omega(0)}{b^*} \rho_d + \f{\tilde \Omega (0)}{\tilde A} \rho_d}.
\end{split}
\end{align}
Relation \eqref{eq:lambda} defines all eigenvalues of the linearized system \eqref{eq:delta:eigen}--\eqref{eq:delta:eigen:3}.
%MD CHANGE 23/04: erase the following paragraph:
%At first sight, the function $f(\lb)$ seems to have three real poles: $\lb=0,$ $\lb=b-b^*$ and
%$$\lb=\alpha' \f{\tilde \Omega(0)}{b^*} \rho_d - \f{\tilde \Omega (0)}{\tilde A} \rho_d =\lb_0.$$
%(Note that $\lb_0 <0$ since we have supposed that exchange functions $\alpha$ and $\omega$ decrease when the population increases.)
%However, taking a closer look at \eqref{eq:lambda} reveals that the only real pole is $\lb=\lb_0$ since the others are cancelled by numerator values, due to the relation between $b$ and $b^*$.
%We also have the following limits:
%$$ \lim_{\lb\to -\infty} f(\lb) = -\infty,\qquad \lim_{\lb\to \lb_0^-} f(\lb) = \pm\infty,
%\qquad \lim_{\lb\to \lb_0^+} f(\lb) = \mp\infty,\qquad \lim_{\lb\to +\infty} f(\lb) = 0,$$
%which show that \eqref{eq:lambda} has at least one real solution.
%Relation \eqref{eq:lambda} that we have obtained here 
It is the same as the characteristic equation obtained from the DDE system \eqref{eq:delay:1}--\eqref{eq:delay:3}. The transition from  stability to instability regions has been studied in many cases of delay equations and we can expect that it corresponds to a Hopf bifurcation that explains the appearance of periodic solutions, see \cite{Diekmann, Magal}. 
Section \ref{section:StabilityRegions} numerically studies this equation, which shows that the nontrivial steady state can either be stable or unstable depending on the parameters of the equation.

%%%%-------------------------------------------------------------------------------------

\section{Stability regions}
\label{section:StabilityRegions}
\setcounter{equation}{0}
\setcounter{figure}{0}
\setcounter{table}{0}

Since the analytic expression of the eigenvalues in Section \ref{section:stabilityanalysisnontrivialstate} is difficult to work out, we determine them numerically. To do so we use the Matlab program DDE-BIFTOOL\cite{engelborghs2001} to determine, equivalently, the eigenvalues of the DDE system \eqref{eq:delay:1}--\eqref{eq:delay:3} as the parameters $\rho_d$ and $b$ vary. In this way, we can determine the regions of stability for the nonzero steady state.  Futhermore, applying Proposition~\ref{prop:Zero}, we can numerically determine where this state disappers and the zero steady state becomes stable.

Notice that this procedure is equivalent to using the eigenvalue equation \eqref{eq:lambda} to numerically determine stability regions for the PDE system (\ref{equation:Approx4alphastar})--(\ref{equation:Approx4omegastarBC}).  However, since the DDE-BIFTOOL software already exists, it is easier to numerically analyze the equivalent DDE system given by \eqref{eq:delay:1}--\eqref{eq:delay:3} instead.

\begin{figure}[htbp]
\begin{center}
\includegraphics[scale = 0.7]{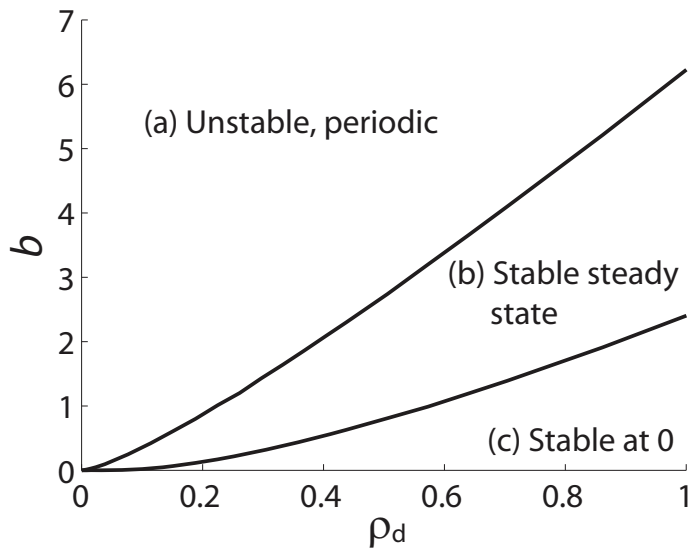}
\end{center}
\caption{
Stability regions for the PDE system from Section~\ref{section:Approximation4} in $(\rho_d, b)$-space.
The variable $\rho_d$ denotes the advection rate of Omega cells with respect to $x$, measured in units of 1/day.  Higher $\rho_d$ means that cells spend less time in the Omega compartment before differentiating.
The variable $b$ is the exponential growth rate of Omega cells, measured in units of 1/day.
(a) In the top region, the system is unstable and exhibits periodic behavior.
(b) In the middle region, the system is stable at the nonzero steady state.  Hence, the Alpha and Omega populations approach a nontrivial equilibrium.
(c) In the bottom region, the nonzero steady state does not exist, and the system is stable at the zero steady state.
}
\label{figure:StabilityRegions}
\end{figure}

Figure~\ref{figure:StabilityRegions} shows the stability regions of the system with respect to $\rho_d$ and $b$.  The variable $\rho_d$ denotes the advection rate of Omega cells with respect to $x$, measured in units of 1/day.  The variable $b$ is the exponential growth rate of Omega cells, measured in units of 1/day.

As we can see from Figure~\ref{figure:StabilityRegions}, the region at which the system is stable at the nontrivial steady state falls between the regions where the system is unstable and where the system is only stable at zero.  Hence, the stability of the system at a nonzero equilibrium depends on a balance between the rate of differentiation and the growth rate of the Omega (proliferating) stem cell population.

If the stem cell growth rate increases or if the rate of differentiation decreases sufficiently, the system transitions to instability, in which case population exhibits periodic behavior.  This transition from a stable equilibrium to unstable periodic behavior might correspond to the transition between CML and acute myelogenous leukemia (AML), a disease characterized by the rapid proliferation and invasion of the blood by immature, undifferentiated cells.  Indeed, unstable and high amplitude oscillations in the stem cell population would result in the sudden overproduction of immature, progenitor cells.  These results concur with those of \cite{adimy2008}, which presents a DDE model of hematopoiesis and concludes that obstructing cell differentation at an early stage of development results in the overproduction of immature cells, potentially leading to AML.

On the other side, if the growth rate decreases or the differentiation rate increases sufficiently, the system loses the nontrivial equilibrium and becomes stable only at zero.  This result makes sense, since the stem cell population has to proliferate fast enough to replenish itself as Omega cells continually differentiate.  However, it is currently unclear to us whether this scenario corresponds to a particular disease.

\begin{figure}[htbp]
\begin{center}
\begin{tabular}{cc}
\includegraphics[scale = 0.55]{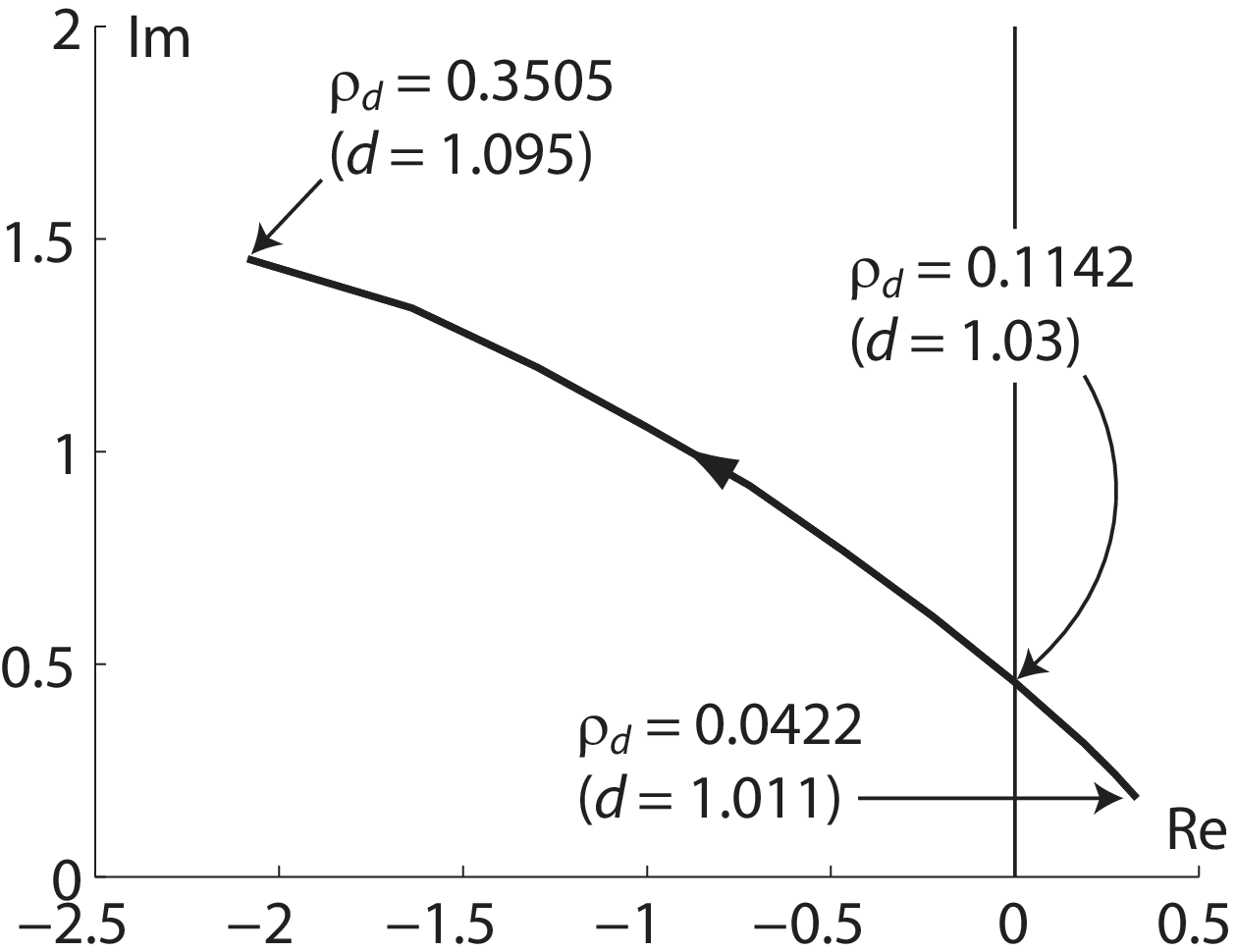} &
    \includegraphics[scale = 0.55]{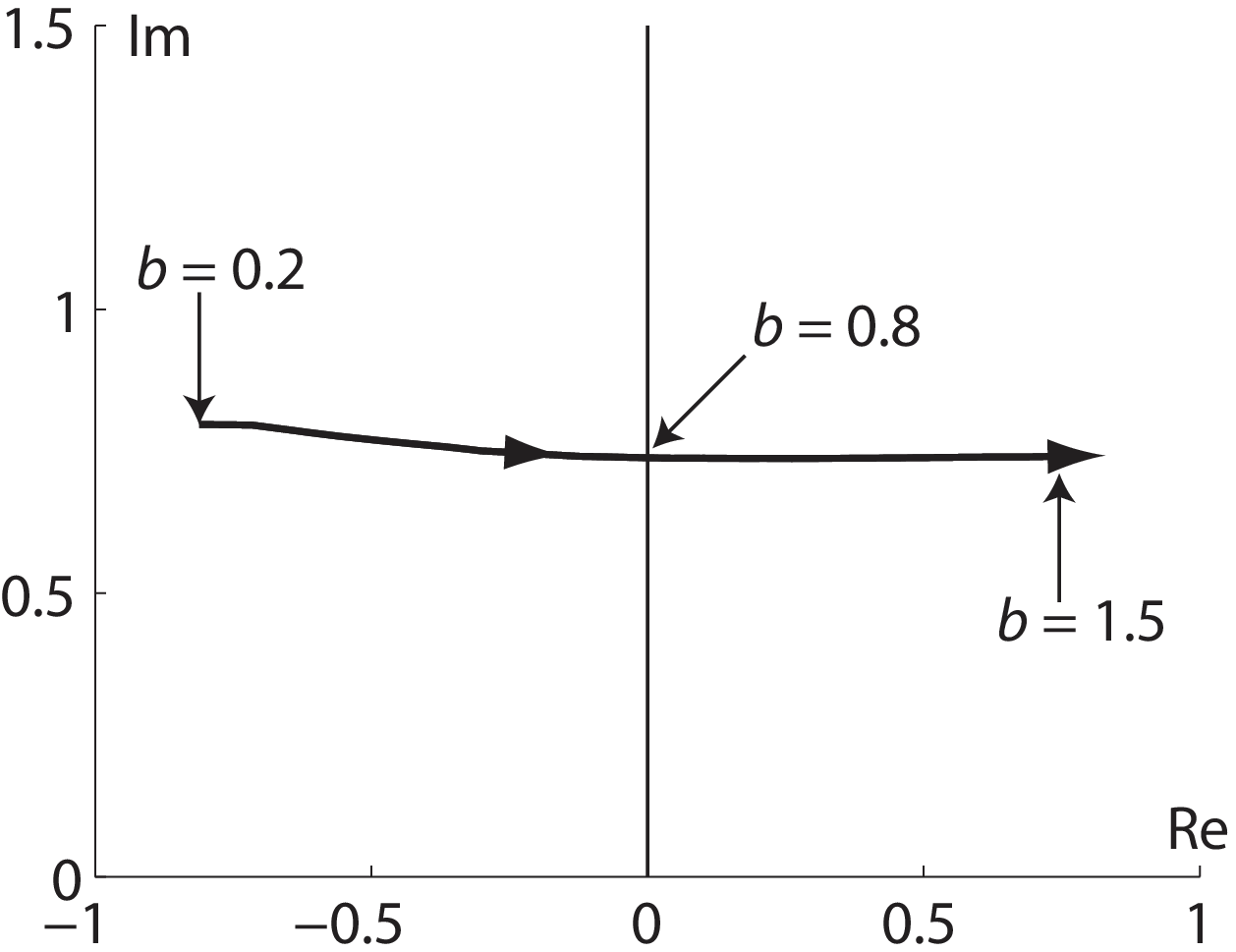} \\
{\bf (a)} & {\bf (b)} \\
\end{tabular}
\end{center}
\caption{
Trajectories of the rightmost eigenvalue around the first stability crossing point.  Each trajectory is associated with a conjugate trajectory (not shown).
(a) Trajectory of the rightmost eigenvalue as $\rho_d$ varies from 0.0422 to 0.3505.  The value of $b$ is fixed at 0.42.
(b) Trajectory of the rightmost eigenvalue as $b$ varies from $0.2$ to $1.5$.  The value of $\rho_d$ is fixed at 0.1884 ($d = 1.05$), which is the estimated value from \cite{roeder2006}
}
\label{figure:RootTrajectories}
\end{figure}

Finally, the program DDE-BIFTOOL can also be used to numerically determine
the trajectories of the rightmost eigenvalue of \eqref{eq:lambda} as it crosses the imaginary axis from stability to instability (i.e., from negative to positive real part) as parameter $\rho_d$ or $b$ is varied.  The trajectories of the rightmost eigenvalues as $\rho_d$ and $b$ are increased independently are shown in Figure~\ref{figure:RootTrajectories}.
As also shown in Figure~\ref{figure:StabilityRegions}, we see that the system moves from instability to stability as $\rho_d$ increases, whereas the system moves from stability to instability as $b$ increases.  

%%%%-----------------------------------------------------------------------------------
%%%%-----------------------------------------------------------------------------------
%%%%-----------------------------------------------------------------------------------
%%%%-----------------------------------------------------------------------------------
%%%%-----------------------------------------------------------------------------------
%%%%-----------------------------------------------------------------------------------
%%%%-----------------------------------------------------------------------------------
%%%%-----------------------------------------------------------------------------------
%%%%-----------------------------------------------------------------------------------
%%%%-----------------------------------------------------------------------------------
%%%%-----------------------------------------------------------------------------------
%%%%-----------------------------------------------------------------------------------
%%%%-----------------------------------------------------------------------------------
\section{Conclusion}
\label{section:Conclusion}
\setcounter{equation}{0}
\setcounter{figure}{0}
\setcounter{table}{0}
In this paper, we simplify the PDE model of hematopoiesis presented in \cite{kimleelevyBMB2008b} and perform a stability analysis.
The PDE model from \cite{kimleelevyBMB2008b} is a time-continuous extension of the ABM of \cite{roeder2006} and has shown to coincide closely with the dynamics of the original ABM.  Our simplifications render the original PDE model amenable to analysis while closely preserving the qualitative behavior.

Prior to our analysis, we simplified the PDE model given by
\eqref{equation:PDEalpha}--\eqref{equation:PDEomegastar} and
\eqref{equation:PDEalphaBC}--\eqref{equation:PDEomegastarBC2}
to the system labelled Approximation 4 in Section~\ref{section:Approximation4}.
As shown in Figure~\ref{figure:ComparisonApprox4toABM}, Approximation 4 is a good approximation of the original ABM from a qualitative perspective.  Furthermore, Approximation 4 can be analyzed analytically.
Hence, the process of simplification and analysis in this paper, allows us to apply an analytical approach to an otherwise inherently complex and intractable ABM.
Furthermore, the observation that Approximation 4 is equivalent to a DDE system shows that the ABM of \cite{roeder2006} possesses inherent similarities to other age-structured DDE models proposed by 
Mackey \emph{et al.} \cite{adimy2005, bernard2003,bernard2004, colijn2005A, colijn2005B, pujo2004} and Adimy \emph{et al.} \cite{adimy2008}.

In addition, Approximation 2, presented in Section~\ref{section:Approximation2}, generates dynamics that are surprisingly close, even quantitatively, to those of the original ABM.  (See Figure~\ref{figure:ComparisonApprox2toABM}.)  Although Approximation 2 is not as suitable for analysis as Approximation 4, it possesses the same computational complexity, since it is also a system of two PDEs with one state variable.  Indeed, it only takes minutes to run numerical simulations of Approximation 2 using the explicit upwind scheme.  This speed of processing is comparable to that of Approximation 4 and of the difference equation model presented in \cite{kimleelevyBMB2008a}.  Hence, for the purposes of numerical simulations, Approximation 2 could be used as an alternative to Approximation 4.  In addition, since Approximation 2 is a PDE system, one could readily add additional ODEs, DDEs, or PDEs to the existing system
to capture the concurrent dynamics of the anti-leukemia immune response during imatinib treatment as was done in \cite{kimleelevyPLoS2008} for the ODE model of imatinib dynamics in \cite{michor2005}.
Since the dynamics of leukemia cells differs between the ABM model of \cite{roeder2006} and the ODE model of \cite{michor2005}, especially in the long-term time scale, it will be useful to examine the impact of the anti-leukemia immune response on the ABM as well.  We leave this for a future work.

Finally, from our analysis, we observe that the system exhibits three types of behavior:
(a) Periodic behavior, (b) A stable nonzero steady state, and (c) Stability at zero.
State (c) corresponds to a state in which the stem cell population cannot maintain itself, and hence dies out to zero.  It would correspond to a condition characterized by insufficient blood production.
State (b) corresponds to the desirable state, in which the hematopoietic stem cells remain at equilbrium and can maintain a stable population level of blood cells.  State (a) corresponds to a state in which stem cell and blood cell populations fluctuate rapidly.  While state (b) captures the behavior of more stable blood cell growth such as in the case of CML, state (a) captures some of the relevant dynamics of AML.  Hence, the transition from state (b) to state (a) might provide insight into the specific malfunction (either in differentiation or growth) that leads to the transition from CML to AML.  
%%%%-----------------------------------------------------------------------------------
%%%%-----------------------------------------------------------------------------------
%%%%-----------------------------------------------------------------------------------
%%%%-----------------------------------------------------------------------------------
%%%-----------------------------------------------------------------------------------
%%%%-----------------------------------------------------------------------------------
\appendix
%%%%-----------------------------------------------------------------------------------
%%%%-----------------------------------------------------------------------------------
%%%%-----------------------------------------------------------------------------------
%%%%-----------------------------------------------------------------------------------
%%%%-----------------------------------------------------------------------------------
\section{Parameter estimates}
\label{appendix:parameterestimates}
\setcounter{equation}{0}
\setcounter{figure}{0}
\setcounter{table}{0}
The sigmoidal transition functions are given in \cite{roeder2006} by
\begin{align}
f_{\alpha/\omega} (\Alphatotal/\Omegatotal) = 24 \left( \frac{1}{\nu_1 + \nu_2 \exp \left( \nu_3 \frac{\Alphatotal/\Omegatotal}{\tilde{N}_{A/\Omega}} \right)} + \nu_4 \right),
\label{equation:sigmoidfunctions}
\end{align}
where $\Alphatotal$ and $\Omegatotal$ denote the total populations in the Alpha and Omega compartments, respectively (see (\ref{equation:PDEtotals}).
Furthermore,
\begin{align*}
\nu_1  = \frac{ h_1 h_3 -h_2^2 }{ h_1 + h_3 - 2h_2}, \qquad
\nu_2  = h_1 - \nu_1, \qquad
\nu_3  = \ln \left(\frac{h_3 - \nu_1}{ \nu_2}\right), \qquad
\nu_4  = f_{\alpha/\omega} (\infty),
\end{align*}
\begin{align*}
h_1 & =  \frac{1}{f_{\alpha/\omega}(0) - f_{\alpha/\omega}(\infty)}, \quad
h_2 & = \frac{1}{ f_{\alpha/\omega}(\tilde{N}_{A/\Omega}/2) - f_{\alpha/\omega}(\infty) }, \quad
h_3 & = \frac{1}{f_{\alpha/\omega}(\tilde{N}_{A/\Omega}) - f_{\alpha/\omega}(\infty) }.
\end{align*}
The values of the various parameters are listed in Table~\ref{table:RoederParameters}.

\vspace{7cm}

\begin{table}[htbp]
\footnotesize
\centering
\begin{tabular}{|l|l|l|lll|l|}
\hline
Param & Description & Ph$^-$ & \multicolumn{3}{|l|}{Ph$^+$/imatinib-affected} & Rescaled \\
\hline
\hline
$a_{\text{min}}$ & Min value of affinity $a$  &  0.002 &  \multicolumn{3}{|l|}{0.002} &   \\
$a_{\text{max}}$ & Max value of affinity $a$  &  1.0   &  \multicolumn{3}{|l|}{1.0}   &   \\
\hline
$d$ & Differentiation coefficient & 1.05 & \multicolumn{3}{|l|}{1.05}  & $\rho_d = 0.1884$ \\
$r$ & Regeneration coefficient    & 1.1  & \multicolumn{3}{|l|}{1.1}   & $\rho_r = 0.3681$ \\
\hline
$c_1$       & Duration of S/G$_2$/M-phases       & 17 hours   & \multicolumn{3}{|l|}{17 hours} & 17/24 \\
$c_2$       & Cell cycle duration                & 49 hours   & \multicolumn{3}{|l|}{49 hours}  & 49/24 \\
\hline
$\lambda_{\text{p}}$  & Lifespan of precursor cells & 20 days & \multicolumn{3}{|l|}{20 days} & 20 \\
$\lambda_{\text{m}}$  & Lifespan of mature cells & 8 days & \multicolumn{3}{|l|}{8 days} & 8 \\
\hline
$\tilde{\tau}_c$  & Division period for precursors   & 24 hours  & \multicolumn{3}{|l|}{24 hours}  & 1 \\
\hline
$f_{\alpha}(0)$             &  Transition characteristic for $f_{\alpha}$ & 0.5    & \multicolumn{3}{|l|}{1.0}   &    \\
$f_{\alpha}(\tilde{N}_A/2)$ &  Transition characteristic for $f_{\alpha}$ & 0.45   & \multicolumn{3}{|l|}{0.9}   &    \\
$f_{\alpha}(\tilde{N}_A)$   &  Transition characteristic for $f_{\alpha}$ & 0.05   & \multicolumn{3}{|l|}{0.058} &    \\
$f_{\alpha}(\infty)$        &  Transition characteristic for $f_{\alpha}$ & 0.0    & \multicolumn{3}{|l|}{0.0}   &    \\
$\tilde{N}_A$               &  Scaling factor & $10^5$ & \multicolumn{3}{|l|}{$10^5$} & 1 \\
\hline
$f_{\omega}(0)$             &  Transition characteristic for $f_{\omega}$ & 0.5    & 1.0  & / & 0.0500  &   \\
$f_{\omega}(\tilde{N}_A/2)$ &  Transition characteristic for $f_{\omega}$ & 0.3    & 0.99 & / & 0.0499  &   \\
$f_{\omega}(\tilde{N}_A)$   &  Transition characteristic for $f_{\omega}$ & 0.1    & 0.98 & / & 0.0498  &   \\
$f_{\omega}(\infty)$        &  Transition characteristic for $f_{\omega}$ & 0.0    & 0.96 & / & 0.0496  &   \\
$\tilde{N}_{\Omega}$               &  Scaling factor & $10^5$ & \multicolumn{3}{|l|}{$10^5$} & 1 \\
\hline
\hline
$r_{\text{inh}}$ & Inhibition intensity  & & \multicolumn{3}{|c|}{0.050} & 1.2 \\
$r_{\text{deg}}$ & Degradation intensity & & \multicolumn{3}{|c|}{0.033} & 0.8 \\
\hline
\hline
$b$ &    Expansion rate for Approximations  & & \multicolumn{3}{|c|}{ } & 0.42 \\
$\kappa$ &  Fraction of time spent in G$_1$-phase  & & \multicolumn{3}{|c|}{ } & 0.54 \\
\hline
\end{tabular}\caption[\hspace{0.05cm} Parameters]{
Parameters from \cite{roeder2006} and rescaled values used in this paper.
The inhibition intensity, $r_{\text{inh}}$, refers to the probability that a proliferative Ph$^+$ cell 
(i.e., an $\Omega$ or $\Omega^*$ cell) becomes imatinib-affected in a given time interval.
The degradation intensity, $r_{\text{deg}}$, refers to the probability that an imatinib-affected, 
proliferative Ph$^+$ cell dies in a given interval.
In the column for rescaled parameters, if the entry is blank, the corresponding parameters are left unchanged.
}
\label{table:RoederParameters}
\end{table}

%%%%-----------------------------------------------------------------------------------

\section{Algorithm for the agent-based model}
\label{appendix:RoederAlgorithm}
\setcounter{equation}{0}\setcounter{figure}{0}
\setcounter{table}{0}
We summarize the algorithm of the ABM from \cite{roeder2006}.  At every time step (1 hour), the ABM is defined as the following set of actions:

\noindent \textbf{A. Preliminary calculations.}
\begin{enumerate}
\item Calculate the total populations of $A$ and $\Omega$ cells.
\item During imatinib treatment:
    \begin{itemize}
    \item Remove the proliferative Ph+ cells ($\Omega^+$ and $\Omega^{+/i}$) that undergo apoptosis.
    \item Determine which unaffected proliferative Ph+, $\Omega^+$, become imatinib-affected.
    \end{itemize}\end{enumerate}

\noindent \textbf{B. Proliferation, death, change of state, clocks.}
At this stage, all cells fall into one of three categories: $A$ stem cells, 
$\Omega$ stem cells, differentiated cells.  \begin{enumerate}
\item For each $A$ stem cell:
    \begin{itemize}    \item Determine whether the cell transfers to $\Omega$.  If a cell transfers, skip the remaining actions for $A$ cells.  
Note that the transition function depends on whether the cell is Ph$^-$, Ph$^+$, or imatinib-affected.  Calculate transition 
probabilities based on the total population of $\Omega$ calculated in Step A1.
    \item Increase the cell's affinity by a factor of $r$.
    \end{itemize}
\item For each $\Omega$ stem cell:    \begin{itemize}
    \item Determine whether the cell transfers to $A$.  If a cell transfers, skip the remaining actions for $\Omega$ cells.  
Calculate transition probabilities based on the total population of $A$ calculated in Step A1.
    \item If the cell's affinity is less than or equal to $a_{\text{min}}$, the cell becomes a differentiated cell of age $0$.  
If the cell differentiates, skip the remaining actions for $\Omega$ cells.
    \item If a cell's affinity is greater than $a_{\text{min}}$, decrease the cell's affinity by a factor of $d$.
    \item Increase the counter $c$ by $1$.
    \item If the counter $c$ is greater than or equal to 49, set $c$ to 0 and create a new cell with identical attributes and 
state values as the current cell.
    \end{itemize}
\item For each differentiated cell:
    \begin{itemize}
    \item Increase the cell's age by one.
    \item If the cell's age is a multiple of $24$ between $24$ and $480$, inclusively, create a new differentiated cell with 
with the same age as the current cell.
    \item If a cells age reaches $672$, that cell dies.
    \end{itemize}
\end{enumerate}
\noindent Note that differentiated cells of age less than 480 are considered to be proliferating precursors, whereas 
differentiated cells of age greater than or equal to 480 are considered to be nonproliferating mature cells.

%%%%-----------------------------------------------------------------------------------

%%%%-----------------------------------------------------------------------------

\bibliographystyle{plain}
\bibliography{references,Bibil_090315}
%%%%-----------------------------------------------------------------------------

\end{document}